\documentclass[final]{elsarticle}       % onecolumn (second format)
\usepackage{graphicx}
\usepackage{latexsym,hyperref,lineno}
\usepackage{lipsum}
\usepackage{amsfonts,amsmath}
\usepackage{hyperref}
\usepackage{graphicx}
\modulolinenumbers[5]
\usepackage{algorithmic}
\usepackage[color]{showkeys}
\definecolor{refkey}{gray}{.75}
\usepackage{booktabs}
\def\p{\partial}
\def\tilde{\widetilde}
\def\hat{\widehat}
\def\bs{\boldsymbol}

\def\ub{\bs{u}}

\def\kb{\bs{\mathrm{k}}}

\def\BC{\mathcal{BC}}

\def\gradh{\nabla_h}
\def\sgradh{\nabla^\perp_h}

\def\M{\mathcal{M}}

\begin{document}

\begin{frontmatter}

\title{Conservative numerical schemes with optimal dispersive wave
  relations -- Part II.~Numerical evaluations
}

%% Group authors per affiliation:
\author{Qingshan Chen, Lili Ju, and Roger Temam}

%% or include affiliations in footnotes:
%\author[mymainaddress,mysecondaryaddress]{Elsevier Inc}

\begin{abstract}
A new energy and enstrophy conserving scheme for the shallow water
equations is evaluated using a 
suite of test cases over the global spherical or bounded domain. The
evaluation is organized around a set of pre-defined 
properties: accuracy of individual opeartors, accuracy of the whole
scheme, conservation of key quantities, control of the divergence variable,
representation of the energy and enstrophy 
spectra, and simulation of nonlinear dynamics. The results confirm
that the scheme is between the first and second order accurate,
and conserves the total energy and potential enstrophy up to the time
truncation errors. The scheme is  capable of producing more
physically realistic energy and enstrophy spectra, indicating that it
can help prevent the unphysical energy cascade towards the 
finest resolvable scales. With an optimal representation of the
dispersive wave relations, the scheme is able to keep the flow close
to being non-divergent, and maintain the
geostrophically balanced structures with large-scale geophysical flows
over long-term simulations. 
\end{abstract}


\end{frontmatter}

%\linenumbers

%\maketitle

\section{Introduction}

%% Pursuit of perfect numerical schemes for large-scale geophysical
%% flows (why)
In long-term simulations of the geophysical flows, where millions of
iterations are involved, numerical accuracy alone cannot ensure an
accurate representation of certain essential statistics about the
flows. Low-order numerical schemes armed with better conservative
properties and dispersive wave representations can sometimes reproduce
more realistic dynamics than higher order numerical schemes do
(\cite{Arakawa1977-og, Arakawa1981-dy}). Many authors have been
inspired to search for numerical schemes that are not necessarily more
accurate, but possess more desirable properties in conservation of
key quantities, such as mass, energy, enstrophy, and/or representation
of the dispersive wave relations
(\cite{Heikes1995-ky,Ringler2002-jx,Chen2013-fa,Ringler2010-sm,Gassmann2012-xb,
  Thuburn2014-ns}). It is in the same spirit that we develop, in Part
I.~of this project (\cite[CJT1 hereafter]{Chen2019-ls}), a new 
numerical scheme that conserves both the total energy and potential
enstrophy and possess the optimal dispersive wave relations on
unstructured meshes over bounded or unbounded domains. This article
constitutes Part II of the project, and its goal is to numerically
evaluate the scheme.

%% Focal questions: advantages of conservation of energy and enstrophy
%% and optimal 
%% representation of dispersive wave relations in the long-term
%% simulations of geophysical flows.
A focal question in the evaluations of the new scheme is what
advantages it has over conventional numerical schemes that lack some
of its properties. More specifically, we would like to understand how
the energy and enstrophy conservations, and the representation of the
dispersive wave relations affect the long-term dynamics
of the simulated flows.
To fully address these important questions requires
extensive 
tests using real-world applications, which are  out of the
reach and scope of the current project, and are left for future
endeavors. In this project, we content ourselves with a preliminary
study of these questions in a somewhat idealized setting. For the
effect of the 
energy and enstrophy conservations, Arakawa and Lamb
(\cite{Arakawa1981-dy}) make use of an idealized test case with a
meridional ridge within a zonal reentrant channel, and shows that, a
few weeks out into the simulations,
their total energy and potential enstrophy conserving scheme, compared with a
non-conserving scheme, is better at eliminating
the spurious inverse energy cascade, and leads to a less noisy and
more organized wind field. Here, we take a similar approach, and
perform, among all the tests, 
a comparative study using an external numerical model, the MPAS-sw
model (for the shallow water equations), and an existing test case
(the shallow water standard 
test case \#5 with mountain topography). On a global sphere, the
MPAS-sw model operates on exactly the same kind of orthogonal meshes
as our model does, and therefore is an ideal choice for
comparison. In this study, we not only examine the vorticity dynamics,
but also compute and compare 
the energy and enstrophy spectra for the two models at different
points during the
simulations.
As for the representation of the dispersive wave
relations, it is commonly accepted it has profound impact on the
geostrophic adjustment process (\cite{Arakawa1977-og,Randall1994-vu}),
but we are not aware of any numerical studies of its effect in long-term
real applications. In this project, we also perform a second comparative
study, comparing our model with the MPAS-sw model in terms of the
control of the divergence variable. The rationale for this study is
simple: a geostrophically balanced flow should stay close to being
non-divergent. 

%% Other more mundane issues explored in this work: numerical
%% accuracy, conservations
With the focus on the potential impacts of the new scheme in real
applications, this numerical study is conducted around the scheme's
properties, rather than around the established test cases, as is
usually done in the literature
(\cite{Chen2013-fa,Ringler2010-sm,Chen2014-cy,Thuburn2015-hi}). Specifically,
this study considers six properties of the numerical scheme, which are
deemed essential for the success of long-term simulations: the
accuracy of the individual operators, the 
accuracy of the whole model, the conservation of key quantities, the
control of the divergence variables, the representation of the energy
and enstrophy spectra, and the simulation of the barotropic
instabilities. Each property is evaluated using well selected or
specially designed test cases. 
%% Layout of the rest of the article
The rest of the article is organized as follows. Section
\ref{sec:scheme} recalls the conservative numerical scheme.
Section \ref{sec:numerics} presents the numerical results. Finally,
some discussions and remarks
are provided in Section \ref{sec:conclu}.

\section{The conservative schemes}\label{sec:scheme}
The vorticity-divergence form of the nonlinear rotating shallow water
equations on a two-dimensional 
manifold $\M$ read
\begin{equation}
  \label{eq:3}
  \left\{
    \begin{aligned}
      &\dfrac{\p}{\p t}\phi + \nabla\cdot(\phi\ub) = 0,\\
      &\dfrac{\p}{\p t}\zeta + \nabla\cdot(q\phi\ub) = 0,\\
      &\dfrac{\p}{\p t}\gamma - \nabla\times(q\phi\ub) = 
      -\Delta\left(g(\phi+b) + K\right).
    \end{aligned}\right.
\end{equation}
Here, $\M$ may be a bounded region on the
global sphere, or the entire global
sphere. The symbol $\phi$ stands for the fluid top-to-bottom
thickness, $\ub$ the 
horizontal velocity field of the fluids, $\zeta \equiv
\nabla\times\ub$ the relative vorticity, 
and $\gamma\equiv\nabla\cdot\ub$ the divergence.

When treated as a single variable (\cite{Salmon2004-tt, Salmon2005-pd,
  Salmon2007-dm}), 
the mass flux $\phi\ub$  has a
Helmholtz decomposition (\cite{Girault1986-sr})
\begin{equation}
  \label{eq:18}
  \phi\ub = \nabla^\perp\psi + \nabla\chi,
\end{equation}
where $\psi$ and $\chi$ are the streamfunction and the velocity potential
respectively, and $\nabla^{\perp} = \kb\times\nabla$ the skewed
gradient operator. In the classical Helmholtz decomposition, $\psi$ is
already assumed to satisfy the homogeneous Dirichlet boundary
condition (\cite{Girault1986-sr}); to enforce  no-flux boundary
condition on the flow, one only need to set the normal
derivative of the  velocity
potential $\chi$ to zero on the boundary, and thus the boundary
conditions on $\psi$ and $\chi$ are
\begin{subequations}\label{eq:19}
  \begin{align}
    \psi &= 0, \qquad {\rm on}\; \p\M,\label{eq:19a}\\
    \dfrac{\p\chi}{\p n} &= 0,  \qquad{\rm on}\; \p\M.\label{eq:19b}
  \end{align}
\end{subequations}
The relation between $(\psi,\,\chi)$ and the vorticity and divergence
$(\zeta,\,\gamma)$ can be easily derived,
\begin{equation}
  \label{eq:24}
  \left\{
    \begin{aligned}
      \nabla\times\left(\phi^{-1}(\nabla^\perp\psi +
        \nabla\chi)\right) &= \zeta,\\
      \nabla\cdot\left(\phi^{-1}(\nabla^\perp\psi +
        \nabla\chi)\right) &= \gamma.
    \end{aligned}\right.
\end{equation}
Under the boundary conditions just given in \eqref{eq:19} and with a
strictly positive 
thickness field $\phi$, the system \eqref{eq:24} is a {\itshape coupled,
self-adjoint,} and 
{\itshape strictly elliptic} system for $(\psi,\,\chi)$.

By substituting the mass flux $\phi\ub$ \eqref{eq:18} into
\eqref{eq:3}, we can 
eliminate the velocity variable from the shallow water system
entirely,
\begin{equation}
  \label{eq:3a}
  \left\{
    \begin{aligned}
      &\dfrac{\p}{\p t}\phi + \Delta\chi = 0,\\
      &\dfrac{\p}{\p t}\zeta + \nabla\cdot(q\nabla^\perp\psi) +
      \nabla\cdot(q\nabla\chi) = 0,\\ 
      &\dfrac{\p}{\p t}\gamma - \nabla\times(q\nabla^\perp\psi) -
      \nabla\times(q\nabla\chi) =  
      -\Delta\left(g(\phi+b) + K\right).
    \end{aligned}\right.
\end{equation}
As it turns out, this is the form that the to-be-derived numerical schemes
bear close resemblance to. 

The inviscid shallow water system \eqref{eq:3} is a Hamiltonian system
(\cite{Salmon1998-eg}), and can be put in the canonical form,
\begin{equation}
  \label{eq:25}
  \dfrac{\p}{\p t} F = \{F,\,H\}. 
\end{equation}
Here, $F$ represents a functional associated with the shallow water
system, and $H$ the Hamiltonian, which is also a functional and is, for
the shallow water system \eqref{eq:3}, given by
\begin{equation}
  \label{eq:20}
   H = \int_\M \left(\dfrac{1}{2} \phi^{-1}\left(
      |\nabla^\perp\psi|^2   + |\nabla\chi|^2 +
      2\nabla^\perp\psi\cdot\nabla\chi\right) 
     + \dfrac{1}{2}g(\phi + b)^2 \right)  d{\mathbf x}. 
 \end{equation}
The Poisson bracket $\{\cdot,\,\cdot\}$ has three
components,
\begin{equation}
  \label{eq:26}
\{F,\,H\} = \{F,\,H\}_{\zeta\zeta} + \{F,\,H\}_{\gamma\gamma} +
\{F,\,H\}_{\phi\zeta\gamma},
\end{equation}
and each component is defined as follows,
\begin{subequations}\label{eq:27}
  \begin{align}
    \{F,\,H\}_{\zeta\zeta} =& \int_\M q J(F_\zeta,\,H_\zeta) d{\mathbf x},\label{eq:27a}\\
    \{F,\,H\}_{\gamma\gamma} =& \int_\M q J(F_\gamma,\,H_\gamma) d{\mathbf x},\label{eq:27b}\\
    \{F,\,H\}_{\phi\zeta\gamma}
     = &\int_\M \left[q(\nabla F_\gamma\cdot\nabla H_\zeta-\nabla
        H_\gamma \cdot\nabla F_\zeta )\right. +\label{eq:27c}\\
         &\qquad\left. (\nabla F_\gamma\cdot\nabla H_\phi-\nabla
       H_\gamma \cdot\nabla F_\phi )\right] d{\mathbf x}.\nonumber
  \end{align}
\end{subequations}
In the above, $F_\zeta$, etc., are  short-hands for the functional
derivatives $\delta F/\delta\zeta$, etc.
The potential vorticity $q$ for the SWEs is given by
\begin{equation}
  \label{eq:27-2}
  q = \dfrac{f+\zeta}{\phi},
\end{equation}
and $J(\cdot,\cdot)$ is the Jacobian operator   defined as
\begin{equation}
  \label{eq:27-1}
  J(a,b) = \nabla^\perp a \cdot\nabla b
\end{equation}
for any two scalar functions $a$ and $b$.
The Jacobian operator is skew-symmetric w.r.t.~its two argument
functions.  

As a consequence of the skew-symmetry of the Jacobian operator and the
permutations present in its third ($_{\phi\zeta\gamma}$) component, 
the Poisson bracket \eqref{eq:26} is skew-symmetric. Therefore, when
$F$ is replaced by $H$ in \eqref{eq:25}, the right-hand side vanishes,
implying that the total energy is conserved. It is also easy to check
that quantities in the form of
\begin{equation*}
  C = \int_\M \phi q^k dx,\qquad\forall\,k\ge 0
\end{equation*}
are singularities of the Poisson bracket, and therefore quantities in
this form, called Casimirs, such as mass, total vorticity, potential
enstrophy, etc., are also conserved in the shallow water
system.

The discrete numerical schemes are derived in two phases: the
Hamiltonian phase and the Poisson phase.
In the Hamiltonian phase, one obtains a discretization of the
Hamiltonian $H$, defined in \eqref{eq:20}, by the
discrete Hamiltonian, still denoted as $H$,
\begin{multline}
  \label{eq:35}
  % H = \int_\M \left\{ \hat\phi_h^{-1}\left|\nabla_h^\perp\tilde\psi_h +
  %   \nabla_h\chi_h\right|^2 + \dfrac{1}{2}g(\phi_h + b_h)^2\right\} d{\mathbf x}.
  H = \int_\M \left\{ \hat\phi_h^{-1}\left( \left|\nabla_h^\perp\psi_h\right|^2 +
    \left|\nabla_h\chi_h\right|^2 +
    \nabla_h^\perp\tilde\psi_h\cdot\nabla\chi_h 
    +\nabla_h^\perp\psi_h\cdot\nabla_h\tilde\chi_h\right) \right. \\
 \left.\hphantom{\nabla_h^\perp\tilde\psi_h\cdot\nabla\chi_h}
  +\dfrac{1}{2}g(\phi_h + b_h)^2\right\} d{\mathbf x}. \quad
\end{multline}
Here, the subscripted variables, such as $\phi_h$, are the discrete
version of the corresponding continuous variables, such as $\phi$. In
CJT1, the following notational convention is used: subscripted
variables (by $_h$, $_i$, 
etc.) are discrete, while un-subscripted variables are continuous;
un-accented discrete variables are defined at cell centers, while
$\hat{\hphantom{\psi}}$ on the top designates edge-defined variables,
and $\tilde{\hphantom{\psi}}$ designates vertex-defined variables. 
The symbols $\nabla_h$ and $\nabla_h^\perp$ represent the finite
difference approximations of the gradient operator $\nabla$ and the
skew-gradient operator $\nabla^\perp$, respectively. The details are
omitted, but can be found in Appendix B of CJT1. 
The factor of $1/2$ in the kinetic energy in \eqref{eq:20} has
disappeared in the discrete version here due to the fact that only one
component is used in the inner product of the vector fields. One also
notes that the skew-symmetry of the Jacobian term is preserved in the
approximation here, which leads to a symmetric elliptic system later. 

In the Poisson phase, one obtains approximations to the Poisson
brackets.
Just like its continuous counterpart, the discrete Poisson bracket
also has three components, 
\begin{equation}
  \label{eq:44}
  \{F,\, H\}_{h} \equiv   \{F,\,
  H\}_{h,\phi\zeta\zeta} +   \{F,\,
  H\}_{h,\phi\gamma\gamma} +   \{F,\,
  H\}_{h,\phi\zeta\gamma}.
\end{equation}
Each component is given by
\begin{align}
  \{F, \,H\}_{h,\zeta\zeta} \equiv{} & 
    \int_\M \hat q_h J_h\left(\dfrac{\delta F}{\delta \zeta_h},\,\dfrac{\delta
  H}{\delta \zeta_h}\right) d{\mathbf x},\label{eq:41}\\
  \{F, \,H\}_{h,\gamma\gamma}
    \equiv{ } &\int_\M \hat 
  q_h J_h\left(\dfrac{\delta F}{\delta \gamma_h},\,\dfrac{\delta
  H}{\delta \gamma_h}\right) d{\mathbf x}, \label{eq:42}\\
  \{F,\,H\}_{h,\phi\zeta\gamma}
  \equiv{} & 
  2\int_\M \hat q_h \left( \nabla_h
        F_{\gamma_h}\cdot\nabla_h H_{\zeta_h} - \nabla_h
        H_{\gamma_h}\cdot \nabla_h
        F_{\zeta_h}\right) d{\mathbf x} + {}\nonumber\\
  &2\int_\M
    \left( \nabla_h
        F_{\gamma_h}\cdot\nabla_h H_{\phi_h} - \nabla_h
        H_{\gamma_h}\cdot \nabla_h
        F_{\phi_h} \right) d{\mathbf x}.& &\label{eq:43}
\end{align}
In the above, where $\hat q_h$ is the potential vorticity at cell
edges, which is a remapping of the PV $q_h$ at cell centers, and $q_h$
is defined (see~\eqref{eq:27-2}) as
\begin{equation}
  \label{eq:27-3}
  q_h = \dfrac{f+\zeta_h}{\phi_h}.
\end{equation}
The discrete geopotential $\Phi_h$, vorticity $\zeta_h$, and
divergence $\gamma_h$ are defined as
\begin{equation}
  \label{eq:36}\left\{
    \begin{aligned}
\Phi_h &=   \widehat{\hat\phi_h^{-2}\left(\left|\nabla_h^\perp 
            \psi_h\right|^2 + \left|\nabla_h\chi_h\right|^2 +
          \nabla^\perp\tilde\psi_h\cdot\nabla\chi_h +  
    \nabla^\perp\psi_h\cdot\nabla_h\tilde\chi_h\right)} + g(\phi_h+
b_h),\\
 \zeta_h &= \nabla_h\times\left(\hat\phi_h^{-1}\sgradh\psi_h\right)
 +
 \dfrac{1}{2}\left(\tilde{\gradh\times\left(\hat\phi_h^{-1}\nabla_h\chi_h\right)}
     +
     \nabla_h\times\left(\hat\phi_h^{-1}\nabla_h\tilde\chi_h\right)\right),\\
 \gamma_h &= 
 \dfrac{1}{2}\left({\gradh\cdot\left(\hat\phi_h^{-1}\nabla_h^\perp\tilde\psi_h\right)}
     +
    \tilde{\nabla_h\cdot\left(\hat\phi_h^{-1}\nabla_h^\perp\psi_h\right)}
  \right) + \nabla_h\cdot\left(\hat\phi_h^{-1}\gradh\chi_h\right).
    \end{aligned}\right.
\end{equation}
The discrete Jacobian operator $J_h$ is a skew-symmetric approximation
to its continuous counterpart,
\begin{equation}
  \label{eq:40}
  J_h(a_h, b_h) \equiv \nabla_h^\perp \tilde a_h\cdot\nabla_h b_h -
  \nabla^\perp_h\tilde b_h \cdot \nabla_h a_h. 
\end{equation}

On a global sphere, no boundary conditions are needed, but
both the second and last equation of \eqref{eq:36} contains
redundancy, and the solutions are not unique, for any solution plus
some constants will still be a solution to the system. To ensure
uniqueness, one can replace one equation from each set with an
equation that sets $\psi_h$ or $\chi_h$ to a fixed value at a certain
grid point. Here, without loss of generality, we pick cell $i=0$,
and set
\begin{equation}
  \label{eq:36a}\left\{
    \begin{aligned}
&\psi_0 = 0,\\
&\chi_0 = 0.
    \end{aligned}\right.
\end{equation}
On a bounded domain, the discrete streamfunction satisfies the
homogeneous Dirichlet boundary conditions, while the homogeneous
Neumann boundary conditions for the velocity potential $\chi_h$ are
only implicitly enforced through the specification of the divergence
operator, and the velocity potential is set to zero at  cell $0$ to
ensure unique solvability of the system. Specifically, the boundary
conditions for the coupled 
elliptic system on a bounded domains are
\begin{equation}
  \label{eq:36aa}\left\{
    \begin{aligned}
&\psi_i = 0,\qquad i\in\BC,\\
&\chi_0 = 0.
    \end{aligned}\right.
\end{equation}

With this Poisson bracket \eqref{eq:44}, a discrete Hamiltonian system
can be constructed,
\begin{equation}
  \label{eq:44a}
  \dfrac{\p F}{\p t} = \left\{F, \, H\right\}_h.
\end{equation}
An energy-conserving (EC) numerical scheme can be obtained by sequentially
setting $F$ in the above to $\phi_i$, $\zeta_i$, and
$\gamma_i$. Details of the derivation can be found in CJT1, and here
we list the EC scheme for the purpose of reference, 
\begin{equation}
  \label{eq:58}
  \left\{
    \begin{aligned}
      \dfrac{d}{dt} \phi_i =
      & -\left[\Delta_h \chi_h\right]_i,\\
      \dfrac{d}{dt}\zeta_i   =
  &-\dfrac{1}{2}\left( \left[\tilde{\nabla_h\cdot\left(\hat q_h 
    \nabla_h^\perp\psi_h\right)}\right]_i +\left[
    \nabla_h\cdot \left(\hat q_h 
    \nabla^\perp_h\tilde\psi_h\right) \right]_i\right) -
\left[\nabla_h\cdot\left(\hat 
    q_h\nabla_h\chi_h\right)\right]_i,\\
\dfrac{d}{dt}\gamma_i =
  &\left[\nabla_h\times\left(\hat
      q_h\nabla_h^\perp\psi_h\right)\right]_i + 
  \dfrac{1}{2} \left( \left[\tilde{\nabla_h\times\left(\hat q_h  
    \nabla_h\chi_h\right)}\right]_i + \left[
    \nabla_h\times \left(\hat q_h 
    \nabla_h\tilde\chi_h\right) \right]_i\right)\\
&- \left[\Delta_h
  \Phi_h\right]_i - \dfrac{1}{4|A_i|}\left(\hat q_{e_1} (\chi_i
  -\chi_{i_1}) + \hat q_{e_2} (\chi_{i_2} - \chi_i)\right).
    \end{aligned}\right.
\end{equation}
The terms preceded by $1/4|A_i|$ in the equation for $\gamma_i$ only
appear for boundary cells, i.e.~$i\in\BC$. 

This scheme just given is guaranteed to conserve the total energy,
thanks to the skew-symmetry of the discrete Poisson bracket
\eqref{eq:44}, but it does not conserve the potential enstrophy. The
reason is this quantity, given as
\begin{equation*}
  Z = \int_\M \phi_h q_h^2 dx,
\end{equation*}
is not a singularity of the discrete Poisson bracket. More
specifically, it fails to nullify the first component
$\{\cdot,\,\cdot\}_{h\zeta\zeta}$ of the Poisson bracket. To remedy
this defect, we follow Salmon (\cite{Salmon2009-xr}), and replace this
component by the trilinear Nambu bracket,
\begin{multline}\label{eq:85}
  \{F,\,H,\,Z\}_{h,\zeta\zeta\zeta} = {}\\
  \dfrac{1}{3}\left( \int_\M \hat
  Z_{\zeta_h}\nabla^\perp_h \tilde{F_{\zeta_h}}\cdot\nabla_h
  H_{\zeta_h} d{\mathbf x} + \int_\M \hat
  H_{\zeta_h}\nabla^\perp_h \tilde{Z_{\zeta_h}}\cdot\nabla_h
  F_{\zeta_h} d{\mathbf x} + \right. \\
  \int_\M \hat
  F_{\zeta_h}\nabla^\perp_h \tilde{H_{\zeta_h}}\cdot\nabla_h
  Z_{\zeta_h} d{\mathbf x} 
   - \int_\M \hat
  Z_{\zeta_h}\nabla^\perp_h \tilde{H_{\zeta_h}}\cdot\nabla_h
  F_{\zeta_h} d{\mathbf x} \\
  \left. {} - \int_\M \hat
  H_{\zeta_h}\nabla^\perp_h \tilde{F_{\zeta_h}}\cdot\nabla_h
  Z_{\zeta_h} d{\mathbf x} - \int_\M \hat
  F_{\zeta_h}\nabla^\perp_h \tilde{Z_{\zeta_h}}\cdot\nabla_h
  H_{\zeta_h} d{\mathbf x} \right).
\end{multline}
This trilinear discrete bracket is skew-symmetric with respect to any
two of its three arguments, and as a consequence, it vanishes when $F$
is set to $Z$. 

The modified Poisson bracket thus has the form
\begin{equation}
  \label{eq:89}
  \{F,\, H\}_h = \{F,\,H,\,Z\}_{h,\zeta\zeta\zeta} + \{F,\,
  H\}_{h,\gamma\gamma} + \{F, \,H\}_{h,\phi\zeta\gamma}. 
\end{equation}
The second and third components remain the same as previously defined
in \eqref{eq:42} and 
\eqref{eq:43}, respectively.
An energy and entrophy conserving (EEC) scheme thus results,
\begin{equation}
  \label{eq:88}
  \left\{
    \begin{aligned}
      \dfrac{d}{dt} \phi_i =
      & -\left[\Delta_h \chi_h\right]_i,\\
      \dfrac{d}{dt}\zeta_i   =
  &   -\dfrac{1}{6}\left[
    \tilde{\nabla_h\cdot\left(\hat q_h \nabla_h^\perp \psi_h -
        \hat\psi_h\nabla_h^\perp q_h\right)} +
    \nabla_h\cdot\left(\hat q_h
      \nabla^\perp_h \tilde\psi_h -\hat\psi_h \nabla_h^\perp \tilde
      q_h  \right)\right]_i - {}\\ 
  &\dfrac{1}{3}\left[\hat{\nabla^\perp_h\tilde\psi_h \cdot\nabla_h q_h
      -\nabla^\perp_h \tilde q_h \cdot \nabla_h\psi_h 
    }\right]_i 
 -
\left[\nabla_h\cdot\left(\hat 
    q_h\nabla_h\chi_h\right)\right]_i,\\
\dfrac{d}{dt}\gamma_i =
& \left[\nabla_h\times\left(\hat q_h\nabla_h^\perp\psi_h\right)\right]_i +
\dfrac{1}{2}\left(\left[\tilde{\nabla_h\times\left(\hat q_h 
    \nabla_h\chi_h\right)}\right]_i + \left[
    \nabla_h\times \left(\hat q_h 
    \nabla_h\tilde\chi_h\right) \right]_i\right)\\
&- \left[\Delta_h
  \Phi_h\right]_i - \dfrac{1}{4|A_i|}\left(\hat q_{e_1} (\chi_i
  -\chi_{i_1}) + \hat q_{e_2} (\chi_{i_2} - \chi_i)\right).
    \end{aligned}\right.
\end{equation}
The term preceded by $1/4|A_i|$ in the equation for $\gamma_i$ only
appears for boundary cells ($i\in\BC$). 

\section{Solution of the elliptic system}\label{s:poisson}
Whether an explicit or implicit time-stepping scheme is used, the
elliptic system 
contained in \eqref{eq:36} will have to be solved at every time
step. Thus, an efficient solution of the elliptic system is vital to
the success of the model. Here, we present an iterative scheme coupled
with the Algebraic Multigrid (AMG) method, which exploit the structure
of the coefficient matrix for the benefit of performance of the
solver. To illustrate the structure, we first write the discrete
elliptic system in the usual matrix-vector form. To avoid introducing
more symbols, we wil call the vector
corresponding to the discrete function $\psi_h$ still by $\psi_h$, and the
vector corresponding to $\phi_h$ still by $\phi_h$. We also denote the
coefficient matrix corresponding to each of the discrete differential
operators each by a letter symbol, as detailed in Table
\ref{tab:symbols}. 

\begin{table}[h]
  \centering
  \begin{tabular}{lll}
  \toprule\noalign{\smallskip}
    Operator & Name & Letter symbol\\
\noalign{\smallskip}\midrule\noalign{\smallskip}
    $\displaystyle\nabla_h$ &  Cell-to-edge gradient & G\\
    $\nabla_h$ & Vertex-to-edge gradient & $\textrm{G}^\#$\\
    $\displaystyle\nabla^\perp_h$ &  Cell-to-edge gradient & S\\
    $\nabla^\perp_h$ & Vertex-to-edge gradient & $\textrm{S}^\#$\\
    $\nabla_h\cdot(\hphantom{x}) $ &  Edge-to-cell divergence & D\\
    $\nabla_h\cdot(\hphantom{x} )$ & Edge-to-vertex divergence & $\textrm{D}^\#$\\
    $\nabla_h\times(\hphantom{x} ) $ &  Edge-to-cell curl & C\\
    $\nabla_h\times(\hphantom{x} )$ & Edge-to-vertex curl & $\textrm{C}^\#$\\
    $\hat\phi_h^{-1}(\hphantom{x} )$  & Edge-to-edge multiplication &
                                                                      H\\
    &&\\
    $\tilde{(\hphantom{x})}$  & Cell-to-vertex mapping & N\\
%%    && \\
    $\tilde{(\hphantom{x})}$ & Vertex-to-cell mapping & M\\
                             & Scaling by cell areas & A\\
  \noalign{\smallskip}\bottomrule\noalign{\smallskip}
  \end{tabular}
  \caption{Discrete differential operators and letter symbols}
  \label{tab:symbols}
\end{table}

Replacing the discrete differential operators in the discrete elliptic
system  (the last two equations of
\eqref{eq:36}) by their corresponding letter symbols, 
we can express the system as
\begin{equation}
  \label{eq:90}
  \left(
    \begin{matrix}
      CHS & \dfrac{1}{2}\left(MC^\# HG + CHG^\#N\right) \\
%      \\
      \dfrac{1}{2}\left(DHS^\# N + MD^\# HS\right) & DHG
    \end{matrix}\right)
  \left(
    \begin{aligned}
      \psi_h\\
      \\
      \chi_h
    \end{aligned}\right) = 
  \left(
    \begin{aligned}
      \zeta_h\\
      \\
      \gamma_h
    \end{aligned}\right)
\end{equation}
This system, which is a discretization of
the elliptic problem \eqref{eq:24}, is not symmetric. It can be made symmetric by a
scaling by the cell areas. Hence, we multiply each row by the diagonal
matrix $A$ representing the cell areas, and we obtain
\begin{equation}
  \label{eq:91}
  \left(
    \begin{matrix}
      ACHS & \dfrac{1}{2}\left(AMC^\# HG + ACHG^\#N\right) \\
%      \\
      \dfrac{1}{2}\left(ADHS^\# N + AMD^\# HS\right) & ADHG
    \end{matrix}\right)
  \left(
    \begin{aligned}
      \psi_h\\
      \\
      \chi_h
    \end{aligned}\right) = 
  \left(
    \begin{aligned}
      A\zeta_h\\
      \\
      A\gamma_h
    \end{aligned}\right)
\end{equation}

Due to the symmetry of the whole system, we have that $ACHS$ and
$ADHG$ are symmetric, and
\begin{equation}
  \label{eq:91a}
(AMC^\#HG 
+ ACHG^\#N)^T = (ADHS^\#N + AMD^\# HS).  
\end{equation}
Due to the fact that $C=D$
and $G=S$, we actually have
\begin{equation}
  \label{eq:92}
  ACHS = ADHG.
\end{equation}
In addition, thanks to the fact that
$C^\# = -D^\#$ and $G^\# = -S^\#$, we have
\begin{align*}
  AMC^\# HG &= - AMD^\#HS,\\
  ACHG^\#N &= - ADHS^\#N.
\end{align*}
Thus, regarding the off-diagonal blocks of \eqref{eq:91}, we have
\begin{equation}
  \label{eq:93}
  AMC^\# HG + ACHG^\#N = - \left(ADHS^\# N + AMD^\# HS\right).
\end{equation}
Combining this relation with \eqref{eq:91a},
we find that
\begin{equation}
  \label{eq:94}
  \left(AMC^\# HG + ACHG^\#N\right)^T = - \left(AMC^\# HG + ACHG^\#N\right),
\end{equation}
that is the off-diagonal block itself is skew-symmetric. 

We denote
\begin{align*}
  P &= ACHS, &  Q & = ADHG,\\
  R & = \dfrac{1}{2}\left(AMC^\# HG + ACHG^\# N\right), & &\\
  X & = \psi_h, &  Y & = \phi_h,\\
  F & = A\zeta_h, &  G & = A\gamma_h.
\end{align*}
Taking \eqref{eq:93} into account, we can write the linear system
\eqref{eq:91} succinctly as
\begin{equation}
  \label{eq:95}
  \left(
    \begin{matrix}
      P & & R\\
      && \\
      R^T & & Q
    \end{matrix}\right) \left(
    \begin{matrix}
      X\\ \\ Y
    \end{matrix}\right) = \left(
    \begin{matrix}
      F\\ \\G
    \end{matrix}\right). 
\end{equation}
This system is a straight transcription of the discrete system
$\eqref{eq:36}_{2,3}$. The boundary conditions, \eqref{eq:36a} or
\eqref{eq:36aa}, can be incorporated into the system by modifying
entries of the 
  coefficient matrices $P$, $Q$, and $R$. Modifying the entries of
  these matrices according to the boundary conditions will retain the
  sysmetry of the linear system \eqref{eq:95}, and render it
  invertible, and thus uniquely solvable. 

For large-scale linear systems, especially those that have to be
solved repeatedly, finding a solver that is efficient in memory and
scale well with the matrix size is essential. Direct solvers, such as
the LU factorization, perform poorly on both fronts, because they have
to store the factorizations, which may be non-sparse, and they
often have a computing load of $O(n^3)$. Krylov-type iterative schemes,
e.g.~CG, are very efficient on memory, but as the matrix size
grows, the convergence slow down dramatically, making them unsuitable
for large-sized ($\ge 40,000$) problems. The Geometric Multigraid Method
(GMG) performs superbly in terms of both memory and speed, if nested
meshes are used and there is an inherent geometric 
structure in the coefficient matrix. Unfortunately, we use generally
unstructured meshes, and such a geometric structure is absent. This
leaves the algebraic multigraid method (AMG) our only choice. The
performance of the AMG is less robost than that of the GMG, but it has
the advantage that it can be used as a black-box solver to any kind of
linear systems.

In our experiments, the AMG is applied to the linear system
\eqref{eq:95}. It is observed that the convergence slows down, as
manifested by the increase in the number of iterations required, as
the matrix size grows. On a small mesh with 10K grid
points (matrix size $\approx$ 20K$\times$ 20K), the AMG requires less
than 100 iterations to converge, while on a large mesh with 600K grid
points 
(matrix size $\approx$ 1.2M$\times$ 1.2M), the AMG requires over 500
iterations. AMG usually performs well on classical problems, such the
Poisson equation. The system \eqref{eq:95} results from the
discretization of  two coupled Poisson equations. Submatrices $P$ and
$Q$ come from discretization of these two Poisson equations, while the
off-diagonal block matrices $R$ is result of the coupling, albeit
weak, between these two. Hence, our hypothesis is that the slow-down
of the AMG is due to the departure of the coefficient matrix from
that of a single Poisson equation, in both the increased size and the
non-zero off-diagonal entries.

To remedy the slow-down of the AMG on the coupled system
\eqref{eq:95}, we introduce an iterative scheme, which works with the
AMG to solve the system. The main idea of this iterative scheme is
that the entries in the off-diagonal matrix $R$ are generally very small
(they would vanish completely if the thikness fied is uniform), and
thus can be moved to the right-hand side, and treated as known
quantities in an 
iterative process. To be precise, we assume that, in the iterative
process, $X^n$ and $Y^n$ are already known. Then an update of these
variables can be found by solving
\begin{subequations}\label{eq:96}
\begin{align}
  &PX^{n+1} = F - RY^n,\label{eq:96a}\\
  &QY^{n+1} = G - R^T X^n. \label{eq:96b}
\end{align}  
\end{subequations}
The AMG method can be employed on these equations separately, which
amounts to solving a single classical Poisson equation. In our
experiments, across a wide range of 
meshes (with grids points $2,000 \sim 600,000$), the scheme
\eqref{eq:96} converges
in about $4\approx 6$ iterations, and within each iteration, the AMG
method takes about 20 iterations. With the iterative scheme
\eqref{eq:96}, we are able to solve the coupled linear system
\eqref{eq:95} with the AMG method with a computational cost that is
approximately linear to the size of the problem.

\section{Numerical evaluations}\label{sec:numerics}
This section present numerical results from the model. Rather than
going through various well-established test cases for the shallow
water model, as has usually been done in the literature, here, we
focus on a few pre-defined properties of numerical models for
large-scale 
geophysical flows, and examine our model with regard to each one of
these 
properties. Specifically, we will examine our model with regard to six
properties: 
(1) the accuracy of individual 
operators in the model; (2) the accuracy of the entire model; (3)
the conservation of some physical quantities; (4) the control of the
divergence 
variable; (5) the representation of the energy and enstrophy spectra;
(6) the simulation of the barotropic instabilities. We believe that
these six properties as a whole are essential for the success of
long-term simulations of large-scale geophysical flows. 

For this numerical study, a selected set of well-established test
cases are used, including the famous shallow water standard test cases
(SWSTC) \#2 (a steady-state zonal flow) and \#5 (zonal flow over a
mountain topography) from \cite{Williamson1992-cq}, and the barotropic
instability test case from \cite{Galewsky2004-hr}. Following
\cite{Bauer2018-uv}, we also use an adapted version of the SWSTC \#2
on the northern hemishpere. Finally, we also make use of a new test
case of our own (\cite{Chen2018-bu}) involving a freely evolving gyre
in the northern Atlantic.

An external model, the MPAS-sw model, is used for
comparison in the examination of certain numerical properties listed
above. The MPAS-sw is a C-grid momentum-based shallow water
model. This model operates on the same kind of unstructured meshes as
our model does (i.e.~spherical centroidal Voronoi tessellations with
Delaunay duals), and has a similar convergence rate (presumed for our
model at 
this point), but it is designed to conserve the total energy only,
not 
the potential enstrophy. Therefore, it is expected both models will
have similar performances in certain aspects, but differ in
others. What these differences mean for long-term simulations of
geophysical flows 
is of course an interesting and important quesiton.

\subsection{Accuracy of individual differential operators}

\begin{table}[h]
% table caption is above the table
\caption{Discrete differential operators. The symbols $\alpha$ and
  $\beta$ stand 
  for generic scalar variables, and $\alpha_h$ and $\beta_h$ generic
  discrete variables defined at cell centers. Other notational
  conventions described in the previous section apply. }
\label{tab:1}       % Give a unique label
% For LaTeX tables use
\begin{tabular}{lllll}
  \toprule\noalign{\smallskip}
  Discrete &&&& Analytic \\
\noalign{\smallskip}\midrule\noalign{\smallskip}
  $\Delta_h\beta_h$
  &  &  & & $\Delta\beta$ \\
\noalign{\smallskip}\hline\noalign{\smallskip}
  $\nabla_h\cdot(\hat\alpha_h \nabla_h\beta_h)$,
  &$\nabla_h\times(\hat\alpha_h \nabla^\perp_h\beta_h)$
  &\hphantom{$\nabla_h\times(\hat\alpha_h \nabla^\perp_h\beta_h)$,}
  &\hphantom{$\nabla_h\times(\hat\alpha_h \nabla^\perp_h\beta_h)$,}
  & $\nabla\cdot(\alpha\nabla\psi)$ \\
  \noalign{\smallskip}\hline\noalign{\smallskip}
  $\tilde{\nabla_h\times(\hat\alpha_h \nabla_h\beta_h)}$,
  &$-\tilde{\nabla_h\cdot(\hat\alpha_h\nabla^\perp_h\beta_h)}$  
  &\hphantom{$\nabla_h\times(\hat\alpha_h \nabla^\perp_h\beta_h)$,}
  &\hphantom{$\nabla_h\times(\hat\alpha_h \nabla^\perp_h\beta_h)$,}
  & $\nabla\times(\alpha\nabla\psi)$ \\
  \noalign{\smallskip}\noalign{\smallskip}
$\tilde{\nabla_h\cdot(\hat\beta_h\nabla^\perp_h\alpha_h)}$
  &\hphantom{$\nabla_h\times(\hat\alpha_h \nabla^\perp_h\beta_h)$,}
  &\hphantom{$\nabla_h\times(\hat\alpha_h \nabla^\perp_h\beta_h)$,}
  &\hphantom{$\nabla_h\times(\hat\alpha_h \nabla^\perp_h\beta_h)$,}
  &\hphantom{ $\nabla\times(\alpha\nabla\psi)$} \\
  \noalign{\smallskip}\hline\noalign{\smallskip}
  $\nabla_h\times(\hat\alpha_h \nabla_h\tilde\beta_h)$,
  &$-\nabla_h\cdot(\hat\alpha_h \nabla^\perp_h\tilde\beta_h)$
  &\hphantom{$\nabla_h\times(\hat\alpha_h \nabla^\perp_h\beta_h)$,}
  &\hphantom{$\nabla_h\times(\hat\alpha_h \nabla^\perp_h\beta_h)$,}
  & $\nabla\times(\alpha\nabla\psi)$ \\
  \noalign{\smallskip}\noalign{\smallskip}
${\nabla_h\cdot(\hat\beta_h\nabla^\perp_h\tilde\alpha_h)}$  
  &\hphantom{$\nabla_h\times(\hat\alpha_h \nabla^\perp_h\beta_h)$,}
  &\hphantom{$\nabla_h\times(\hat\alpha_h \nabla^\perp_h\beta_h)$,}
  &\hphantom{$\nabla_h\times(\hat\alpha_h \nabla^\perp_h\beta_h)$,}
  &\hphantom{ $\nabla\times(\alpha\nabla\psi)$} \\
  \noalign{\smallskip}\noalign{\smallskip}
  $-2\hat{\nabla^\perp_h\tilde\beta_h \cdot\nabla_h \alpha_h}$
  &\hphantom{$\nabla_h\times(\hat\alpha_h \nabla^\perp_h\beta_h)$,}
  &\hphantom{$\nabla_h\times(\hat\alpha_h \nabla^\perp_h\beta_h)$,}
  &\hphantom{$\nabla_h\times(\hat\alpha_h \nabla^\perp_h\beta_h)$,}
  &\hphantom{ $\nabla\times(\alpha\nabla\psi)$} \\
  \noalign{\smallskip}\noalign{\smallskip}
$2\hat{\nabla^\perp_h\tilde\alpha_h \cdot\nabla_h \beta_h}$
  &\hphantom{$\nabla_h\times(\hat\alpha_h \nabla^\perp_h\beta_h)$,}
  &\hphantom{$\nabla_h\times(\hat\alpha_h \nabla^\perp_h\beta_h)$,}
  &\hphantom{$\nabla_h\times(\hat\alpha_h \nabla^\perp_h\beta_h)$,}
  &\hphantom{ $\nabla\times(\alpha\nabla\psi)$} \\
  \noalign{\smallskip}\bottomrule\noalign{\smallskip}
\end{tabular}
\end{table}

\begin{figure}[h]
  \centering
  \includegraphics[width=5in]{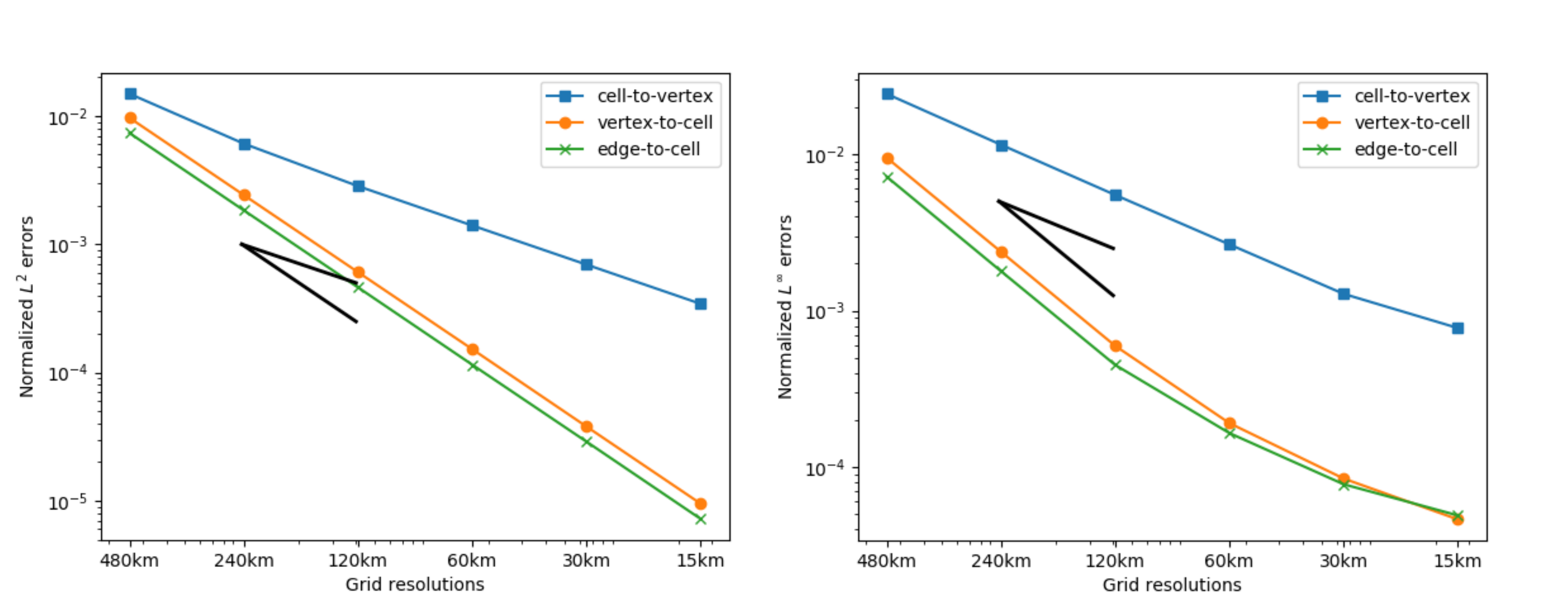}
  \caption{Accuracies of the remapping operators as measured in the
    $L^2$ (left) and $L^\infty$ (right) norms. The reference first
    and second-order convergence curves (black solid lines) are also
    shown. }
  \label{fig:map}
\end{figure}

List in Table \ref{tab:1} are all the discrete differential operators
that appear in either the energy-conserving scheme \eqref{eq:58} or
the energy and enstrophy-conserving scheme \eqref{eq:88}.  The
analytical differential operators, of which these discrete
differential operators are, formally,
approximations of, are also listed in the rightmost column of the
table. This table contains four groups of differential operators,
as divided by the solid lines. 
Discrete differential operators on the same row are
algebraically, as well as numerically, equivalent, while numerical
evidence shows that the operators in the same group are very close in
accuracy (see below). There are a total
of eight discrete differential operators approximating the analytical
differential operator $\nabla\times(\alpha\nabla\psi)$. These eight
operators are divided into two groups. The first group
of three operators (the third group from the table) all have the
vertex-to-cell mapping as the last 
operation, while the other group (the fourth group from the table) has
the same mapping but as the very 
first operation. 

All the operators given in Table \ref{tab:1}, except the discrete
Laplace operator $\Delta_h$, are combinations of three to four
individual discrete operators, e.g.~the discrete gradient, the
discrete divergence, the cell-to-vertex mapping, and the
vertex-to-cell mapping.
% The 
% accuracy of each individual operator is easy to obtain, but the
% analysis of the accuracy of the combinations of three to four of these
% operators is still absent.
% A rigorous analysis of the accuracies of these operators, especially
% those of the compositions of these individual operators, is still
% absent. 
It is reasonable to expect that these
composite operators will remain second-order accurate on a planar
centroidal Voronoi mesh with perfect hexagons. However, on a regular
centroidal 
Voronoi mesh which consists of mostly hexagons and a few pentagons,
the accuracy will likely be 
degraded. Rigorous numerical analysis of these operators on
unstructured meshes will be carried
out elsewhere. Here, we examine the accuracy of these operators
numerically. 

Due to the prominent roles that the remapping (cell-to-vertex,
vertex-to-cell, and edge-to-cell) operators play throughout the
schemes, we examine the accuracies of these operators first. The
cell-to-edge mapping is not included in this study, because the
accuracy of this mapping is well understood to be
second-order accurate, due to the fact that every two neighboring cell
centers are equi-distant to the common Voronoi edge that separates
them. We measure the accuracy of the remapping operatros by comparing
the remapped values from one location to another with the actual
values of a pre-chosen scalar analytic function. Specifically, we use
a scalar  analytic function with mild variations in both latitude and
longitude, 
\begin{equation*}
  \psi = a\cos^3\theta \sin(4\lambda),
\end{equation*}
where $\theta$ and $\lambda$ stand for the latitude and longitude,
respectively, and $a = 6371000\textrm{m}$ the earth radius.  Then
accuracy of the cell-to-vertex remapping operator, for example, is
then computed as follows,
\begin{equation*}
  \dfrac{\|\mathcal{R}_{h,v}\psi - \tilde{\mathcal{R}_{h,c}\psi}\|}{
    \|\mathcal{R}_{h,v}\psi \|}.
\end{equation*}
In the above,
$\mathcal{R}_{h,c}$, and $\mathcal{R}_{h,v}$ stand for the restriction
operators at the cell centers and cell vertices, respectively, and
$\tilde{\hphantom{\psi}}$ the cell-to-vertex remapping, as in the
specification of the schemes.

For this test, we use a series of quasi-uniform spherical centroidal
Voronoi tessellations (SCVT) with resolutions ranging from $480$km up
to $15$km. 
The $L^2$ norm and the $L^\infty$ norm of the errors from all three
remapping operators are calculated, and plotted in Figure
\ref{fig:map}. The $L^2$ norms (left panel) of the errors in the
edge-to-cell and vertex-to-cell remapping operators all converge
towards zero consistently at the second order, while those of the
cell-to-vertex operator converge only at the first order. This
degradation in the accuracy of the cell-to-vertex remapping operator
is likely due to the fact that it is a upscaling remapping from cells
($\times 1$) to vertices ($\times 2$).  The $L^\infty$  errors in all
three remapping operators start out similarly as the $L^2$ errors, up
to the resolution of $120$km. But, onto higher resolutions, the
 curves for the vertex-to-cell and edge-to-cell operators
start to tilt up, and approach the first-order convergence rate after
the $30$km resolution. The curve for the cell-to-vertex operator,
which is already at the a lower convergence rate, exhibits a slight
degradation at the highest resolutions. The degradation in the
accuracies of all three operators are likely indicators of
irregularities in the mesh at the higher resolutions.

\begin{figure}[h]
  \centering
  \includegraphics[width=5in]{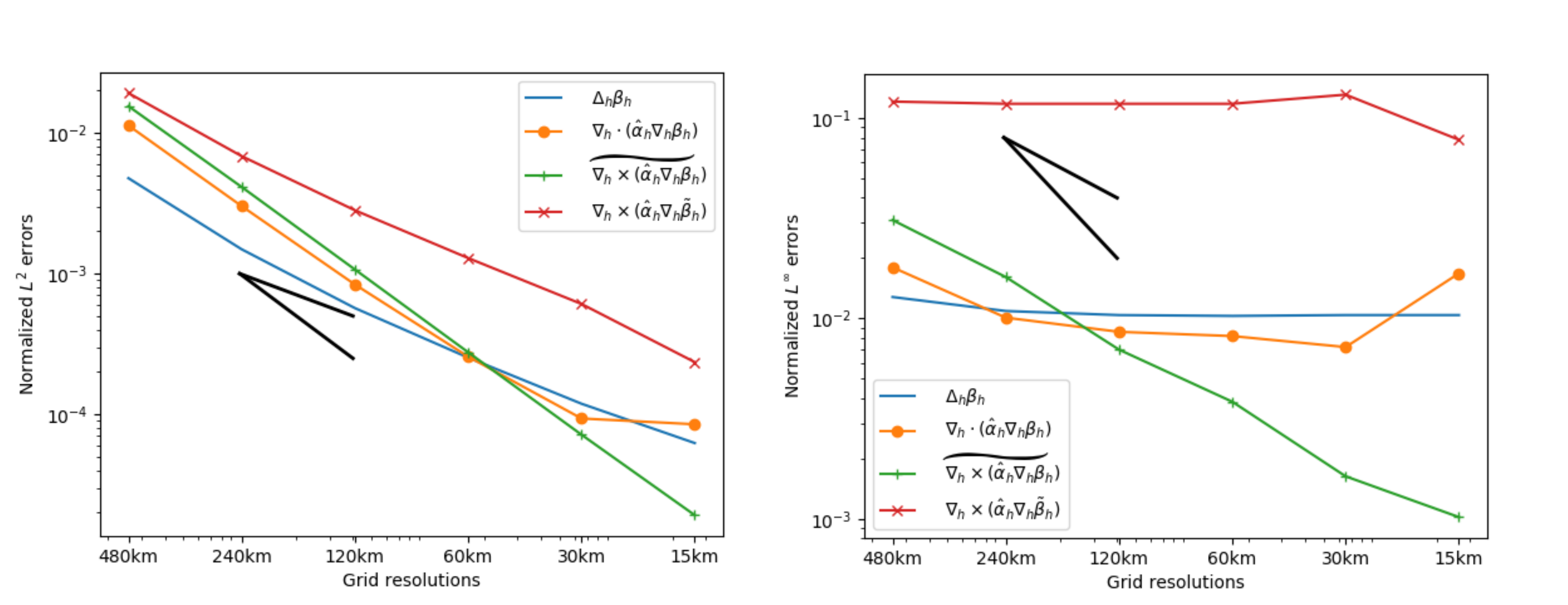}
  \caption{The normalized errors of individual discrete operators,
    along with the first-order and second-order reference convergence
    curves (solid black lines). The errors in 
    the $L^2$-norm of all four representative operators converge at a
    rate between the first and second ordres, while the errors in the
    $L^\infty$-norm do not converge, except that of the representative
    operator from the third group,
%  $tilde{\nabla_h\times(hat\alpha_h \nabla_h\beta_h)}$,
 which converges
    at approximately the first order. }
  \label{fig:operator}
\end{figure}

We study the accuracies of these discrete differential operators of
Table \ref{tab:1} by
measuring their normalized truncation errors. To do
so,  
we take  
\begin{align*}
  &\alpha = 10^{-8}\sin\theta\cos\lambda,\\
  &\beta = a\cos^3\theta\sin 3\lambda,
\end{align*}
where, again, $\theta$ and $\lambda$ stand for the latitude and
longitude, 
respectively, and $a$ the earth radius. The magnitudes of these
variables are chosen so that the magnitude of $\alpha$ matches the
planetary potential vorticity, while that of $\psi$ matches
streamfunction. The 
discrete variables $\alpha_h$ and $\beta_h$ are defined by applying
the restriction operator $\mathcal{R}_h$ to the analytic variable 
$\alpha$ and $\beta$, respectively. The normalized truncation error of
$\nabla_h\cdot(\hat\alpha_h \nabla_h\beta_h)$, for example, is defined
as
\begin{equation*}
  \dfrac{\| \nabla_h\cdot(\hat\alpha_h \nabla_h\beta_h) -
    \mathcal{R}_h\left(\nabla\cdot(\alpha\nabla\beta)\right)\|}
  {\|\mathcal{R}_h\left(\nabla\cdot(\alpha\nabla\beta)\right)\|}.
\end{equation*}

The accuracies of all the discrete operators in Table \ref{tab:1} are
studied, using a series of mesh resolutions ranging from
$480\textrm{km}$ up to $15\textrm{km}$. Numerical results show that
the normalized truncation 
errors for operators in the same group are very close, agreeing up to
at least the third digits. For this reason, we only plot the results
from the first representative member in each group in Figure
\ref{fig:operator}.
On the left are the normalized $L^2$-errors, and
it is seen that the errors from all four representative operators
converge at a rate between the first and the second orders, with the
representative operator
$\tilde{\nabla_h\times(\hat\alpha_h\nabla_h\beta_h)}$ from the third
  group being the most accurate (at the 2nd order) and the
  representative operator
  ${\nabla_h\times(\hat\alpha_h\nabla_h\tilde\beta_h)}$ from the fourth
  the least accurate (at the 1st order).
  The performances of the operators from the first and second groups
  are in between these two extremes.
  The difference in the performances of the third and fourth groups
  are caused by the applications of the remapping operators.
  The operators in the fourth group all start with the cell-to-vertex
  mapping, while the operators in the third group all finish with the
  vertex-to-cell mapping. Both mappings are area weighted, but, as
  shown in Figure \ref{fig:map}, the
  vertex-to-cell mapping is more regular, as it is a downscaling
  mapping from vertices ($\times 2$) to cells ($\times 1$), while the
  cell-to-vertex 
  mapping is more singular, as it is a upscaling mapping from cells
  ($\times 1$) to verteices ($\times 2$). 
    The convergences for the operators from all groups are
    consistent across the whole range of mesh resolutions, except the
    operator $\nabla_h\cdot(\hat\alpha_h\nabla_h \beta_h)$ at the
    highest resolution 
    $15\textrm{km}$. We suspect that this flat-out is caused by an
    anomaly in 
    the mesh. 

On the right of Figure \ref{fig:operator} are the normalized
$L^\infty$ errors. The $L^\infty$ errors for most of the operators are
not converging, which is expected. It is therefore surprising to see
that the $L^\infty$ errors for the representative operator
$\tilde{\nabla_h\times(\hat\alpha_h\nabla_h\beta_h)}$ from the
third group actually converges at the first order. This operator also
has the best and most consistent $L^2$ converging rate. 

The above results concerning the individual operators of Table
\ref{tab:1}, except the third group, largely agree with those of
\cite[Section 
6.2]{Eldred2015-yg} on the SCVT grids in that the operators are
between first and second order accurate in the $L^2$ norm, and
non-consistent in the $L^\infty$ norm. In addition to the SCVT, other
types of optimized grids are also considered in the  work just referenced,
and one particular type, the so-called tweaked iscosahedra, is shown
to lead to consistent operators in both the $L^2$ and $L^\infty$
norms. Operators in the third group on Table \ref{tab:1} are not
needed, and hence not studied in \cite{Eldred2015-yg}. 

\subsection{Accuracy of the entire model}
\subsubsection{A stationary zonal flow on a global sphere (SWSTC\#2) and
  on a hemisphere}
The SWSTC \#2 specifies an initial zonal flow that is in
perfect geostrophic balance. Thus, this test case not only provides a
rare means of testing the accuracy of numerical scheme with an
analytic solution to the full SWEs, it is also an excellent test on
the capability of a numerical scheme to maintain the 
geostrophic balance. In order to test the impact of the presence of a
boundary on the accuracy of the numerical scheme as well as its
capability to maintain geostrophic balance, we follow
\cite{Bauer2018-uv} and run the same
test, with the same initial conditions, but on the northern hemisphere,
with the equator serving as the southern boundary.  
Both the EC scheme \eqref{eq:58} and the EEC scheme \eqref{eq:88} are
used for this test case, and 
the results are very similar. Hence we only present and discuss the
results from the EEC scheme.

The exact configurations of the SWSTC \#2 can be found in the original
Williamson et al paper (\cite{Williamson1992-cq}) and many other
numerical papers that follow. We run the test for 5 days, as
recommended in the original paper, on a suite
of quasi-uniform meshes on the global sphere, with resolutions ranging
from $480$km up to $30$km. Then we compare the
solutions on day 5 
with the intial state, and compute the normalized $L^2$ and
$L^\infty$ errors. The same test is then carried out again, on a suite
of quasi-uniform mesh on the northern hemisphere. 

\begin{figure}[h]
  \centering
  \includegraphics[width=5in]{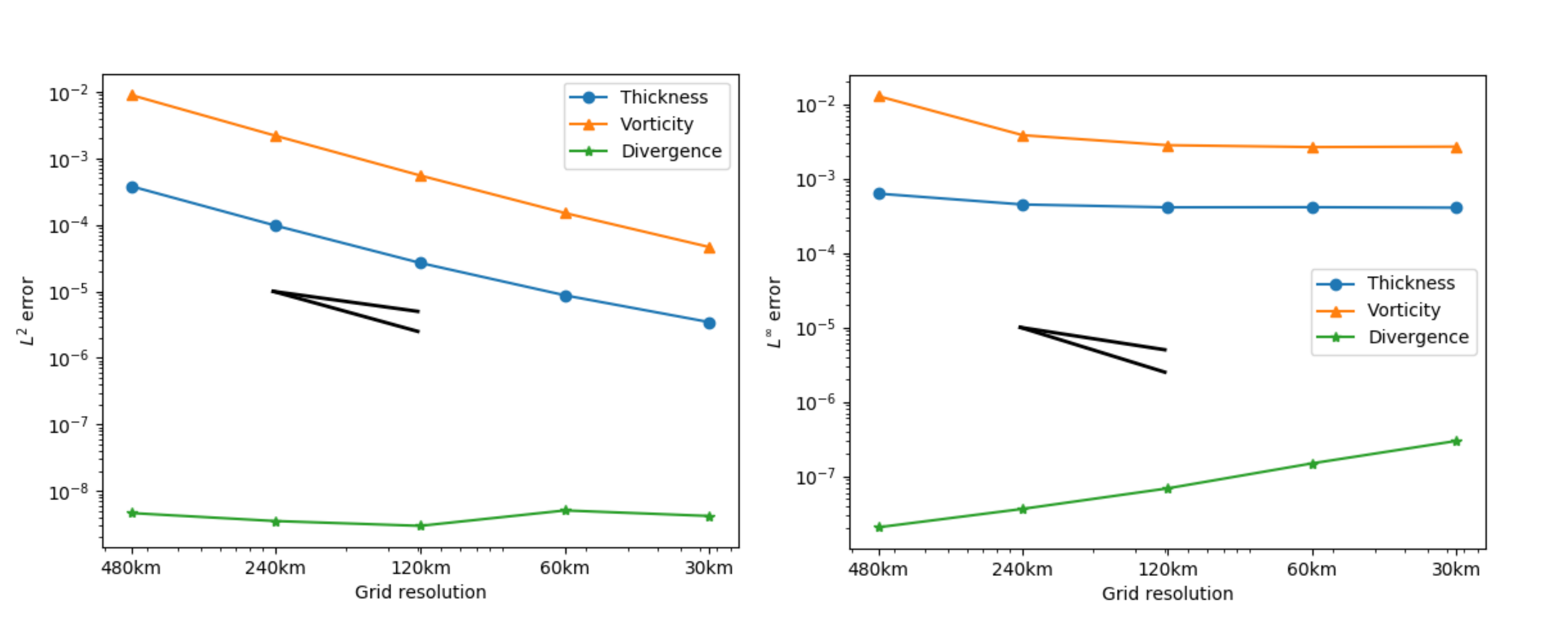}
  \caption{SWSTC \#2: Convergence in the $L^2$ (left) and $L^\infty$
    (right) norms. The errors in the thickness and vorticity variables
  are relative, while the errors in the divergence variables are
  the absolute errors, and their $L^2$ norms are area-weighted and
  -normalized. }  
  \label{fig:swstc2-conv}
\end{figure}

The results on the global sphere are plotted against the grid
resolutions in Figure
\ref{fig:swstc2-conv}. 
% The errors are plotted against the grid resolutions, as represented by the
% numbers of grid cells, in Figure \ref{fig:swstc2-conv},
together with
the first- and second-order reference convergence curves. It is seen
that the $L^2$ errors for both the thickness and vorticity variables
converge at the the second order, with the vorticity errors larger but
converging more consistently, all the way up to the 30km
resolution. Since the
flow starts out as non-divergent, and hence the area-weighted and
-normalized  absolute
errors  for the
divergence variables are used. The $L^2$ errors for
the divergence variable stay below $10^{-8}$ for
all resolutions. The relative $L^\infty$ errors of neither the
thickness or the 
vorticity show any sustained converging trend, and the $L^\infty$
error of the divergence variables shows a steady growing trend. At
30km, the maximum magnitude of the divergence variables reaches about
$10^{-7}$.  

For comparision, the thickness variable under the extended Z-grid
scheme of Eldred 
\cite{Eldred2015-yg} on this particular test case converges at
approximatley the first order in 
both the $L^2$ and $L^\infty$ norms, on the so-called tweaked
icosahedral grids. Under our scheme, the same variable converges
slightly faster in the $L^2$ 
norm on general SCVT grids, but does not converge in the $L^\infty$
norm, consistent with the convergence behaviors of the individual
operators. The convergence behavior of the vorticity variable is not
discussed 
in the work just referenced.

\begin{figure}[h]
  \centering
  \includegraphics[width=5in]{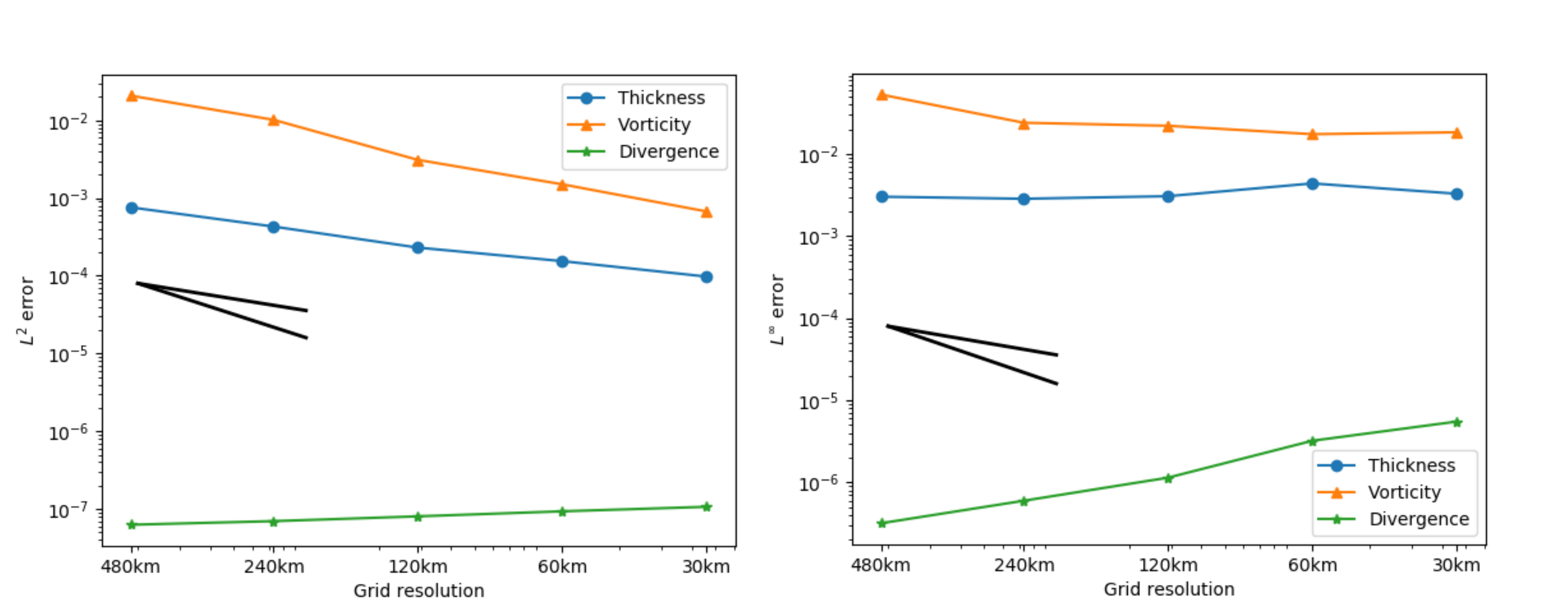}
  \caption{A stationary flow on the northern hemishere: convergence in
    the $L^2$ (left) and $L^\infty$ 
    (right) norms. The errors in the thickness and vorticity variables
  are relative, while the errors in the divergence variables are
  area-weighted and -normalized absolute errors.} 
  \label{fig:swstc12-conv}
\end{figure}

The result from the hemisphere case are plotted
against the grid resolutions in Figure \ref{fig:swstc12-conv},
together with the first- and 
second-order convergence reference curves. In this case, the relative
$L^2$ errors 
in the thickness converge at an approximately first-order rate.
The relative $L^2$ errors for the vorticity variable are larger, but
converge faster, at an approximately 
second-order rate. The $L^\infty$ errors in neither variable converge
consistently. As we already pointed out, the flow starts out
non-divergent, and the magnitude of the divergence variable can only
go up. The figure shows that magnitude in the divergence variable also
go up as the grid refines, with the area-weighted and -normalized
$L^2$ errors reaching at approximately $10^{-7}$, and the $L^\infty$
reaching $5\times 10^{-6}$ at the 30km resolution. 

Compared with the SWSTC \#2 on the global sphere, the errors in the
current case are larger by a factor of 2 to 5, across all
variables and for both norms, and the convergences are also
slower. The downgraded accuracies can be attributed to the presence of
boundaries and the downgraded
uniformity in the grid resolution. The spherical centroidal Voronoi
tessellations on the global sphere are generated using the isohedra as
starting points, and the resulting meshes are highly uniform, with a
ratio of $1.26$ between the highest and lowest resolutions, as
represented by the cell-to-cell distances. The spherical centroidal
Voronoi 
tessellations on bounded domains start with a set of prescribled
points on the boundary and a set of largely random points in the
interior. With iterations, the meshes can achieve a ratio of $2.54$
between the highest and lowest resolutions, which qualifies as
quasi-uniform meshes, but certainly worse than the meshes on a global
sphere.

\subsubsection{Zonal flow over a mountain topography (SWSTC \#5)}
This classical test case starts with an initially zonal flow similar
to that of the SWSTC \#2. The flow impinges on a mountain
topography centered at latitude $30^\circ$ and longitude
$-90^\circ$, and gradually evolves into a turbulent flow around day
25. The exact configurations of this test case, again, can be 
found in the reference \cite{Williamson1992-cq}, and will not be
repeated here. This test case was originally used by Takacs
(\cite{Takacs1988-tx}) to study the conservation of integral
invariants by {\it a posteriori} methods. Subsequently, this
test case is often used by authors to qualitatively study the
capability of numerical schemes to simulate the nonlinear dynamics in
geophysical flows (\cite{Ringler2010-sm, Weller2012-md, Chen2013-fa,
  Li2015-nb}).

\begin{figure}[h]
  \centering
  \includegraphics[width=4in]{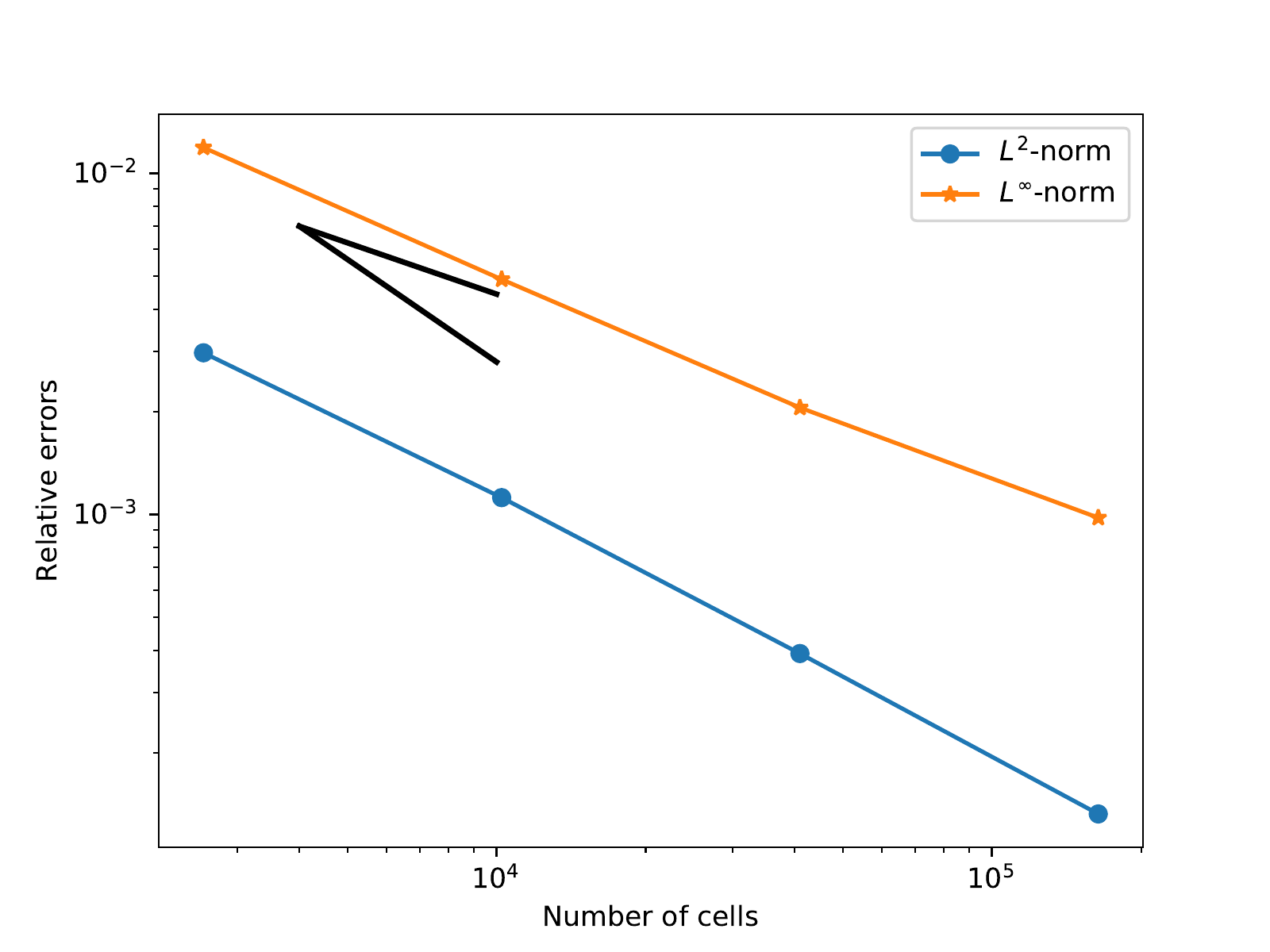}
  \caption{SWSTC \#5: Errors on day 15 with respect to high-resolution spectral
    (H213) solutions. The EEC scheme is used. The grid resolutions
    range from $480$km up to $60$km.}
  \label{fig:swstc5-conv}
\end{figure}

% In this study, we will examine both the conservation and the
% dynamics issues, by comparing  our numerical scheme with the MPAS-sw
% model. But first, we will 
% perform an accuracy check on the thickness variable using a high
% resolution 
% spectral solution (\cite{Jakob-Chien1995-ci}). 
% We compute the errors in the thickness field by comparing the
% solutions to the high-resolution (H213) spectral reference solution
% (\cite{Jakob-Chien1995-ci}).
For this study, we again use a suite of quasi-uniform meshes over the
global 
sphere, with resolutions ranging from $480$km up to $30$km. The
relative errors
in the thickness field on day 15, under both the $L^2$ and the
$L^\infty$ norms, are computed 
by comparing our solutions to the high resolution spectral solution
H213 from (\cite{Jakob-Chien1995-ci}). 
The relative errors are then plotted
against the grid resolutions in Figure \ref{fig:swstc5-conv}, along
with the 
reference first- and second-order convergence curves. It is seen that
the errros in both norms converge consistently at a rate between the
first and second orders. It should be 
pointed out that the reference solution is computed with diffusions,
while our solutions are computed without any 
diffusions. 

\subsection{Conservative properties}
It has been proven in CJT1 that 
the numerical scheme conserves the first order moments, such as the
mass 
and the absolute vorticity, up to the machine errors, and the total
energy and enstrophy, which are second-order moments, up to the time
truncation errors. These results are confirmed in the numerical
study. Here we present the results from two tests. The first test
is on a bounded domain, and the results from this test illustrate the
impact of the time truncation errors in the conservation of the total
energy and the potential enstrophy.
The second test is on a global sphere, and involves comparison to a
third-party model that is designed to conserve the total energy but
not the potential enstrophy. 

\begin{figure}[!ht]
  \centering
  \includegraphics[width=3.2in]{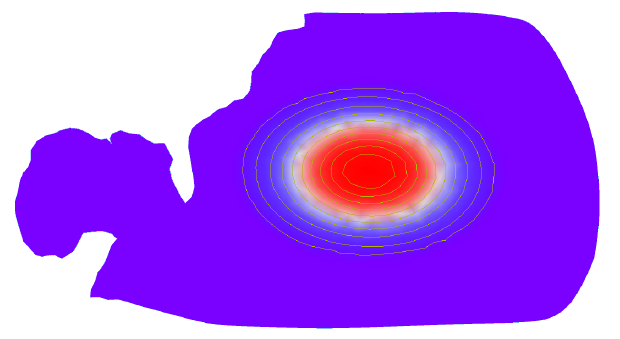}
  \caption{The initial stream function for a test of the pure PV
    advection scheme.}
  \label{fig:trspt-test}
\end{figure}

In numerical schemes, the boundary is often a source of
errors. Inconsistent treatment of the boundary conditions can even
result in unstable simulations. In order to numerically
verify the conservative properties of our numerical scheme and examine
its treatment of the boundary conditions, we design a test case with a
freely evolving gyre in a bounded domain, with no external forcing or
diffusion. Physically, both the total energy and the potential
enstrophy should be conserved, alongside the mass and total
vorticity, and this should also be the case in numerical simulations
with schemes that are designed to conserve these quantities. For the
physical domain, we consider one section of the mid-latitude northern
Atlantic ocean.  
The initial state of the flow is given by 
\begin{equation*}
  \psi(\lambda,\theta) = e^{-d^2} \times (1- \tanh(20*(d-1.5)),
\end{equation*}
where
\begin{equation*}
   d = \sqrt{\dfrac{(\theta -
      \theta_c)^2}{\Delta\theta^2} + \dfrac{(\lambda -
      \lambda_c)^2}{\Delta\lambda^2}},
\end{equation*}
and the parameters $\theta_c = 0.5088$ and $\lambda_c = -1.1$ are the
latitude and longitude of the center point of the Atlantic section,
$\Delta\theta = .08688$, $\Delta\lambda = .15$. The $\tanh$ function
is used to ensure that the initial stream function is flat (with one
constant 
value) near and along the boundary of the domain. A plot of the domain
and 
the stream function is shown in Figure \ref{fig:trspt-test}. This
stream function produces a circular clockwise velocity field in the
middle of the ocean, with the
maximum speed of about $.8$m/s near the center, which is considered
fast for oceanic flows. No external forcing or diffusion is applied,
and the flow is allowed to evolve freely. The circular pattern will
breakdown eventually, because of the irregular shape of the domain and
because of the non-uniform Coriolis force.

\begin{figure}[h]
  \centering
  \includegraphics[width=5in]{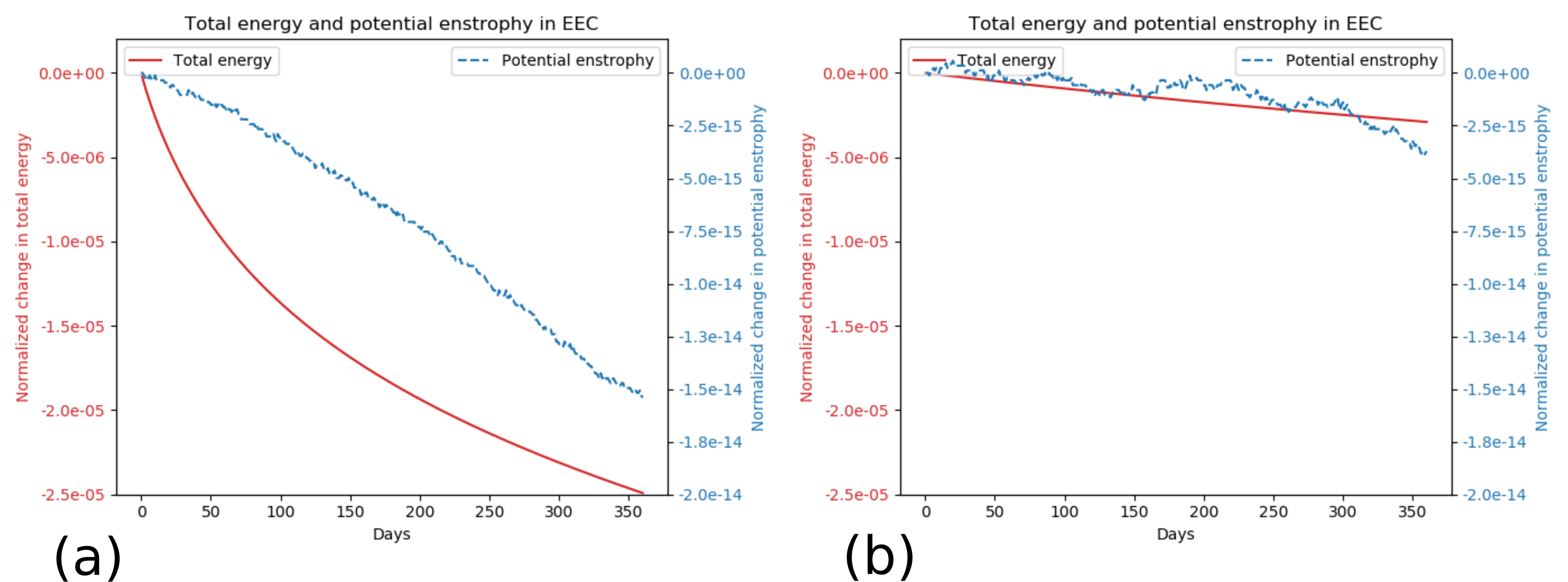}
  \caption{Normalized changes in the total energy and potential
    enstrophy in a freely evolving gyre: (a) with a 3088-cell SCVT
    grid and time step size 80s, and (b) with a 3088-cell SCVT grid
    and time step size 40s. The losses in the total energy and
    potential enstrophy come down with reduced time step sizes,
    confirming that the conservations are up to the time truncation
    errors.}
  \label{fig:free-gyre}
\end{figure}

We cover this bounded domain with a non-uniform 3088-cell SCVT mesh
with 
resolutions ranging from about $20$km to $160$km. The EEC
 scheme is employed, together with the 4th order Runge-Kutta
time stepping scheme, and the simulation is run for one calendar year
(360 days). The total energy and potential enstrophy are computed for
each day, and these statistics are then plotted against time in
Figure \ref{fig:free-gyre}. In order to illustrate the impact of the
time truncation errors, we run the simulation twice, with two time
step sizes, $80$s (left) and $40$s (right). At the end (day 360) of
the first run with 
the $80$s time step size, the system loses $2.5\times 10^{-5}$ of the
total energy and $1.5\times 10^{-14}$ of the potential enstrophy
(Figure \ref{fig:free-gyre} (a)), with the latter almost in the range
of 
roundoff errors for a simulation of this length. 
At the end of the  second run with the $40$s time step size, the
system loses 
$3\times 10^{-6}$ of the total energy, a reduction compared to the
first run by a factor of about 8, and it loses about $4\times
10^{-15}$ of the potential enstrophy, a reduction by a factor of about 4.
These result confirm that the EEC numerical scheme conserve the total
energy and potential enstrophy up to the time truncation errors, and
with reduced time step sizes, the truncation errors do come down.

\begin{figure}[h]
  \centering
  \includegraphics[width=5in]{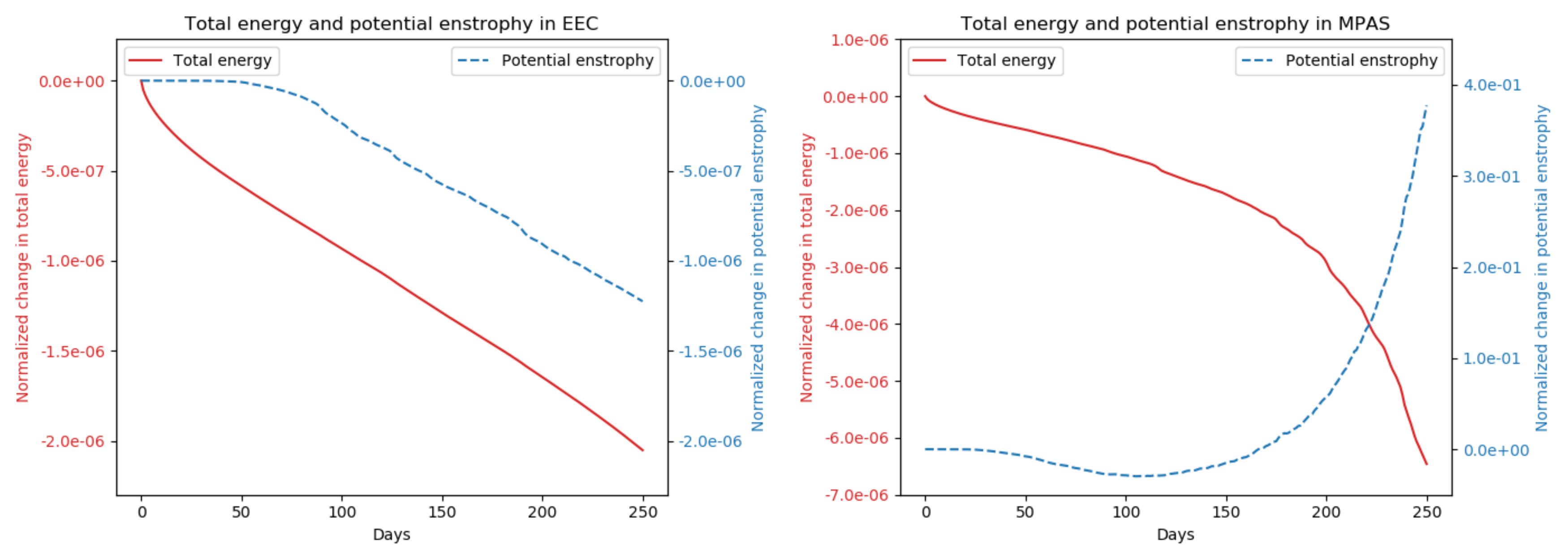}
  \caption{SWSTC \#5: The normalized changes in the total energy and
    the potential enstrophy in EEC (left) and MPAS (right) models.}
  \label{fig:tc5-conserv-EPE}
\end{figure}

Our second test involves an initially zonal flow over a mountain
topography over the global sphere(SWSTC \#5,
\cite{Williamson1992-cq}). The purposes of this test are two
folds. First, it seeks to verify the conservative properties of the
EEC scheme in different setting (flow over topography in a global
domain), and second, it seeks to contrast the results from our EEC
scheme with those of the MPAS-sw model, which operates on the same
kind of  global SCVT meshes, has similar convergence rates,  but is not
designed to conserve the potential enstrophy. For this test, both
models use 
a global quasi-uniform SCVT
mesh with 40962 cells (res.~120km) and the RK4 time stepping
scheme. The simulations are run for $250$ days, and the total energy
and potential enstrophy are computed for each day, and then plotted
against time in Figure \ref{fig:tc5-conserv-EPE}. 
In the EEC scheme (left panel) the flow experiences a weak
decay in both 
the total energy and the potential enstrophy, likely due to the
numerical diffusion of the Runge-Kutter 4th-order time stepping
scheme. By the end of the simulation, the flow loses about  $2\times
10^{-6}$ of the total energy and about $1.2\times 10^{-6}$ of the
potential enstrophy.
In the MPAS-sw model, the evolution curve of the total energy is
similar to that of our EEC scheme, but the potential enstrophy
increases by about 40\% 
at the end of the simulation. This shows that numerical accuracy alone
cannot guarantee the conservation of all dynamical quantities.
Our EEC scheme and the
MPAS-sw model share a similar convergence rate, but have dramatically
different behaviors in terms of potential enstrophy. 

It is worth noting that, our EC scheme, while not proven to conserve
the potential 
enstrophy, manages to only accumulate about $0.8\%$ of it (results not
plotted here),
highlighting the promising prospect of the Hamiltonian approach for
numerical scheme design.

\subsection{Control of the divergence variable}
\begin{figure}[h]
  \centering
  \includegraphics[width=5in]{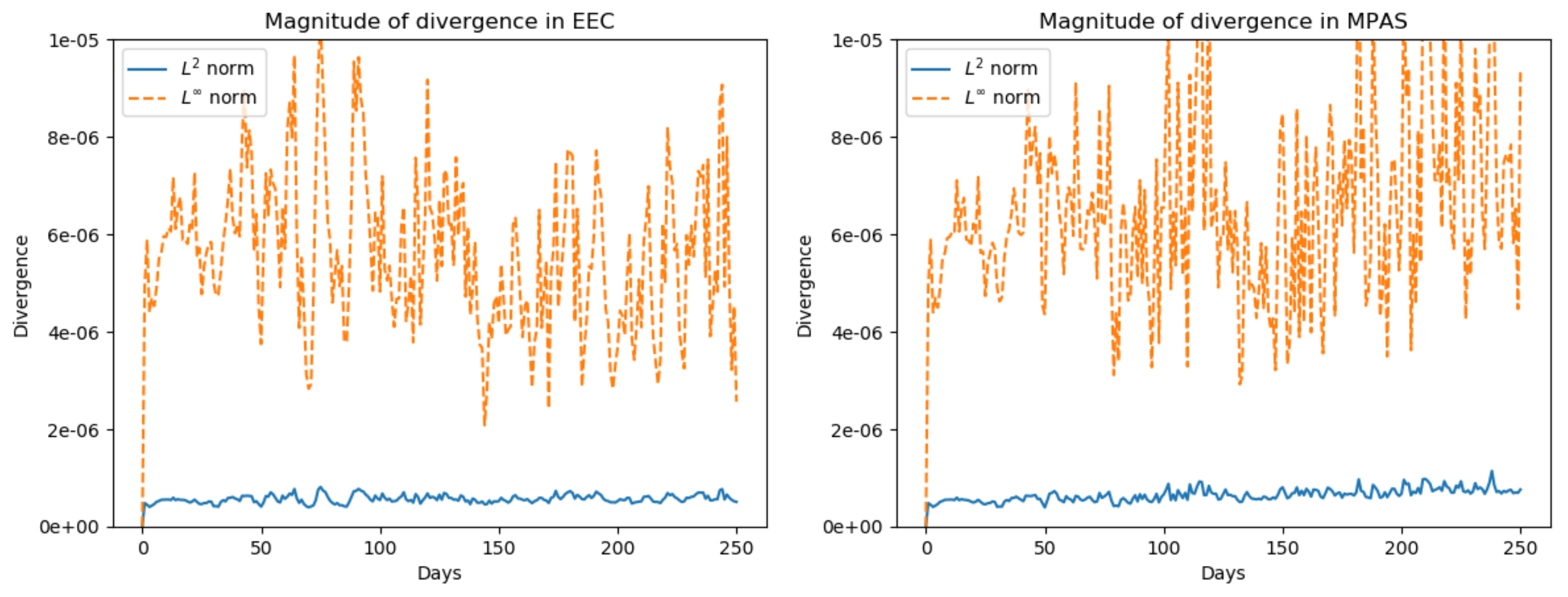}
  \caption{SWSTC \#5: Evolutions in the magnitude in the divergence
    from the EEC (left) and MPAS (right) models.}
  \label{fig:tc5-conserv-div}
\end{figure}

Large-scale geophysical flows constantly evolve around an approximate
geostrophical balance, i.e.~a balance between the Coriolis force and
the pressure gradient. Therefore, it is crucial for numerical scheme
to maintain this geostrophic balance in order to be able to accurately
simulate the dynamics. It has long been recognized that an accurate
representation of the dispersive wave relations is crucial to the
success in this task, because a proper representation of the
dispersive wave relations allow energies to propagate properly across
scales rather unphysically accumulate at one scale, and eventually
destroy geostrophic balance (\cite{Arakawa1977-og,
  Randall1994-vu}). Besides the excellent conservative properties, our
EEC scheme also possess the optimal dispersive wave relations. How
does this last property contributes to its capability to simulate and
maintain geostrophically balanced flows? We need an observable
indicator. When
flows evolve 
around the geostrophic 
balance, the flow is close to being, but not exactly,
non-divergent. Therefore, the divergence variable can be an indicator
to the capability of a numerical scheme to maintain the geostrophic
balance. Ideally, the divergence variable should remain small
throughout the 
simulation. In the left panel of Figure \ref{fig:tc5-conserv-div}, we plot
the evolution of the area-normalized $L^2$ norm and the $L^\infty$ norm of
divergence variable, again on the $120$km quasi-uniform global mesh,
up to day 250. It is seen that the $L^2$ norm of the divergence
variable remains largely flat, exhibiting no noticeable growth. The
$L^\infty$ norm is larger, as expected, but remains largely below
$1\times 10^{-5}$. Again, no noticeable and sustained growth is
observed. For comparison, we also plot the results by the MPAS-sw
model on the same scale, on the right panel of the figure. MPAS-sw is
a C-grid scheme, with decent but 
less accurate (vs 
Z-grid) representations of  the
dispersive wave relations. In the plot, the $L^2$ norm of the
divergence exhibits a slow but steady growth. The $L^\infty$ norm of
the divergence frequently exceeds the upper limit $1\times 10^{-5}$ of
the panel, and, despite being highly oscillatory, also exhibits a
noticeable growth over time. 

\subsection{Representation of the energy and enstrophy spectra}
\begin{figure}[h]
  \centering
  \includegraphics[width=4.6in]{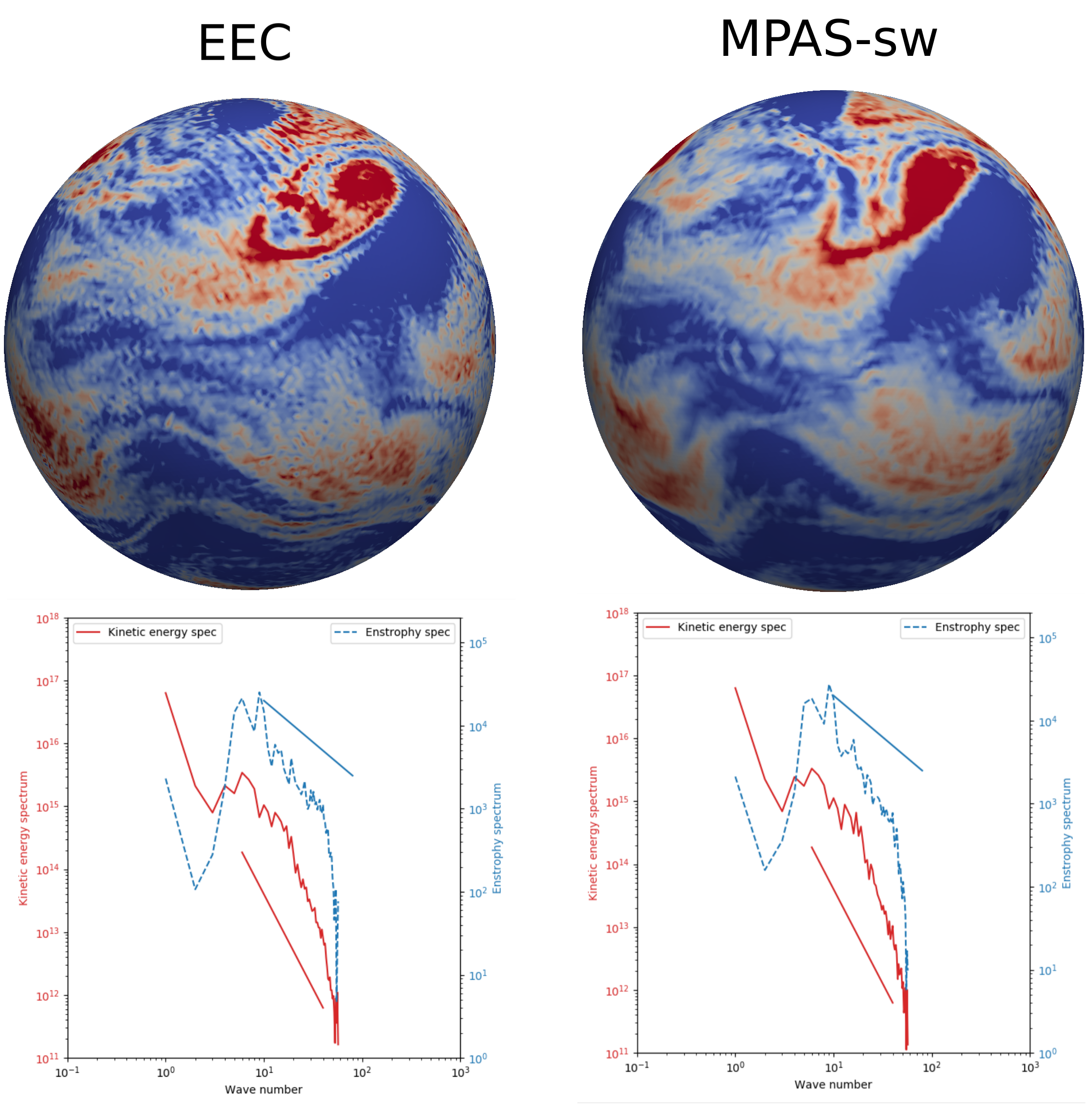}
  \caption{Vorticity field and the energy and enstrophy spectra on day
    50.}
  \label{fig:spectra-day50}
\end{figure}

Our EEC scheme excels in combining the desirable conservative and
dispersive properties into one scheme. Dynamically, what real
advantages does this scheme have compared to other numerical schemes
that are currently in operation? Answers to this big question can
only be sought by applying the scheme to real-world
applications. These important questions will be pursued in the
future. Here 
we examine one crucial dynamical aspect of the numerical scheme, namely its
representation of the spectra of the kinetic energy (KE) and potential
enstrophy (PEns), and the results are compared with those from the
MPAS-sw model, to see what advantages, if any, the EEC model
have. According to the phenomenological theories of 
two-dimensional turbulence (\cite{Batchelor1969-fo, Kraichnan1967-un,
  Lilly1969-wl}), the KE spectrum of a turbulent 
two-dimensional has a slope of $-3$ in the inertial range, while the
spectrum of the PEns has a slope of $-1$. These theories and
predictions have been largely confirmed by numerical studies; see
e.g.~\cite{Chen2011-jh, Maltrud1991-wa}. While a correct
representation of the KE and PEns spectra does not guarantee a
faithful simulation of the nonlinear dynamics of turbulent flows, it
is a prerequisite.

\begin{figure}[h]
  \centering
  \includegraphics[width=4.6in]{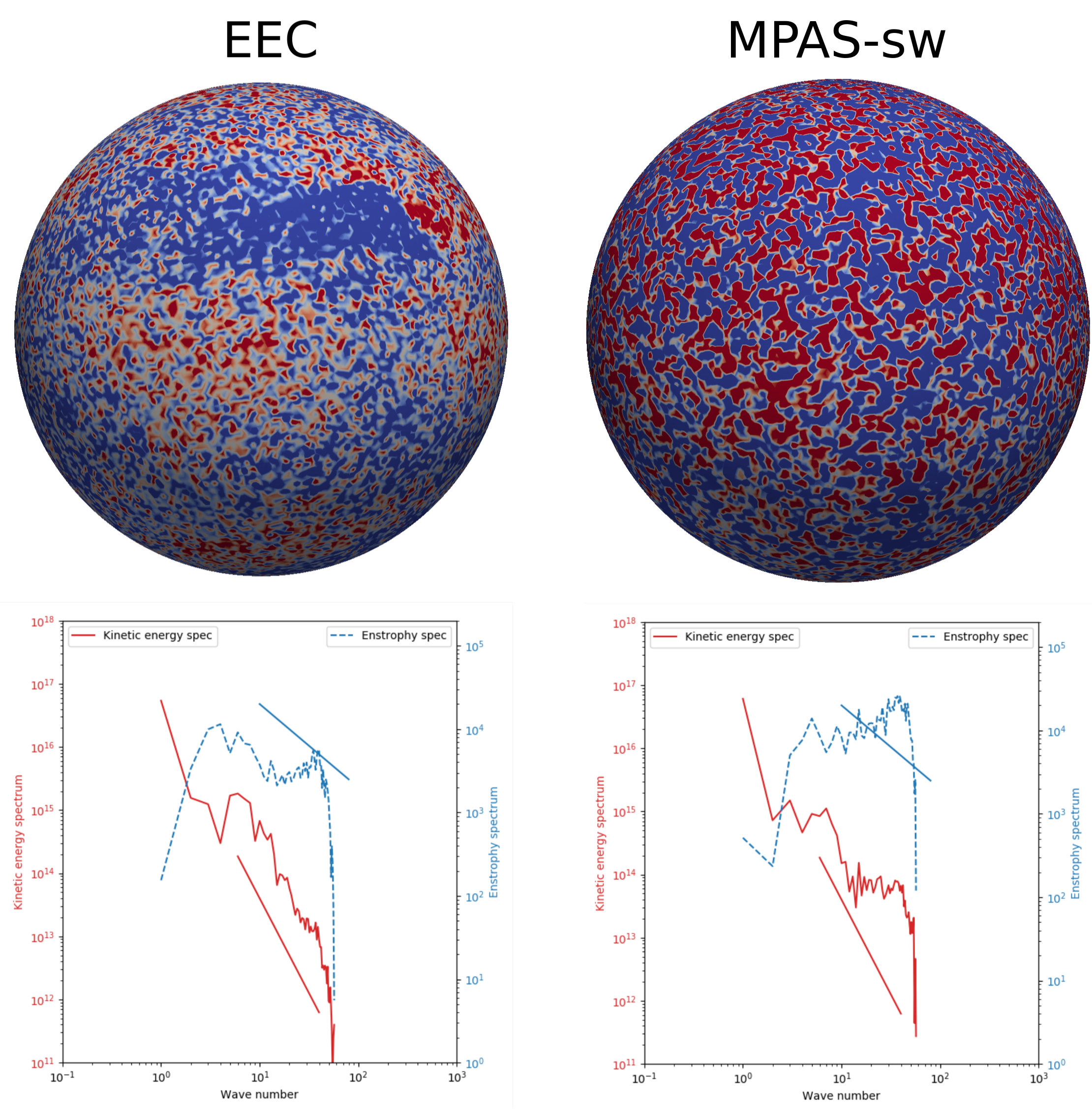}
  \caption{Vorticity field and the energy and enstrophy spectra on day
    250.}
  \label{fig:spectra-day250}
\end{figure}

For this study, we again use the SWSTC \#5 with a free zonal flow over
a mountain topography. We run both EEC and MPAS-sw over the same
global quasi-uniform $120$km mesh for 250 days, and we examine and
plot the vorticity field and the KE and PEns spectra on day 50 and day
250. 
In Figure \ref{fig:spectra-day50} are snapshots of the vorticity field
on day 50 from the EEC and the MPAS-sw models. Placed below the
snapshots are the plots of the KE and PEns spectra at the same
instantaneous moment. The vorticity fields from these two models are
qualitatively similar, while the vorticity field of the MPAS-sw
appears smoother, due to the fact that, in MPAS-sw, the vorticity is
natively defined at the vertices, and is mapped to the cell centers
for the visualization purpose. For the EEC model, the KE spectrum
demonstrates a slope of $-3$ from wave number 10 to 20, while the PEns
spectrum demonstrates a slope of $-1$ from wave number 10 to 30,
approximately verifying the phenomenological theories. The KE and PEns
from the MPAS-sw model exhibit similar trends. For both models, the KE
and PEns spectra decay quickly after wave number 30, apparently due to
the fact that fluid motions still largely concentrate at large scales,
and ``cascade'' is not complete yet. This completely changes on day
250 (Figure \ref{fig:spectra-day250}). For the EEC model, the KE
spectrum now display a slope of $-3$ 
from wave number 10 up to wave number 40, indicating that a
significant portion of the fluid motions has shifted from the large
scales to small scales. Due to this shifting, and due to the lack of
diffusion, the PEns spectrum is now distorted with the PEns piling up
at high wave numbers, and with no discernible part of the curve
conforming to the $-1$ slope. The results from the MPAS-sw are
worse. With this model, the KE spectrum is completely distorted with no 
discernible part conforming to the $-3$ slope, and so is the PEns
spectrum. For the PEns spectrum, not only is the slope $-1$ is gone,
the peak density at high wave number is much larger than the peak
density of the PEns from the EEC model, which is most likely due to
the spurious source of PEns that we saw in the right panel of Figure
\ref{fig:tc5-conserv-EPE}. Snapshots of the vorticity fields
corroberate what the KE and PEns spectra tell us. Both snapshots
are noisy, but the one from the EEC model clearly still possess some
large-scale structures, while the one from the MPAS-sw is completely
noisy, with no discernible structures. 

\begin{figure}[h]
  \centering
  \includegraphics[width=3.7in]{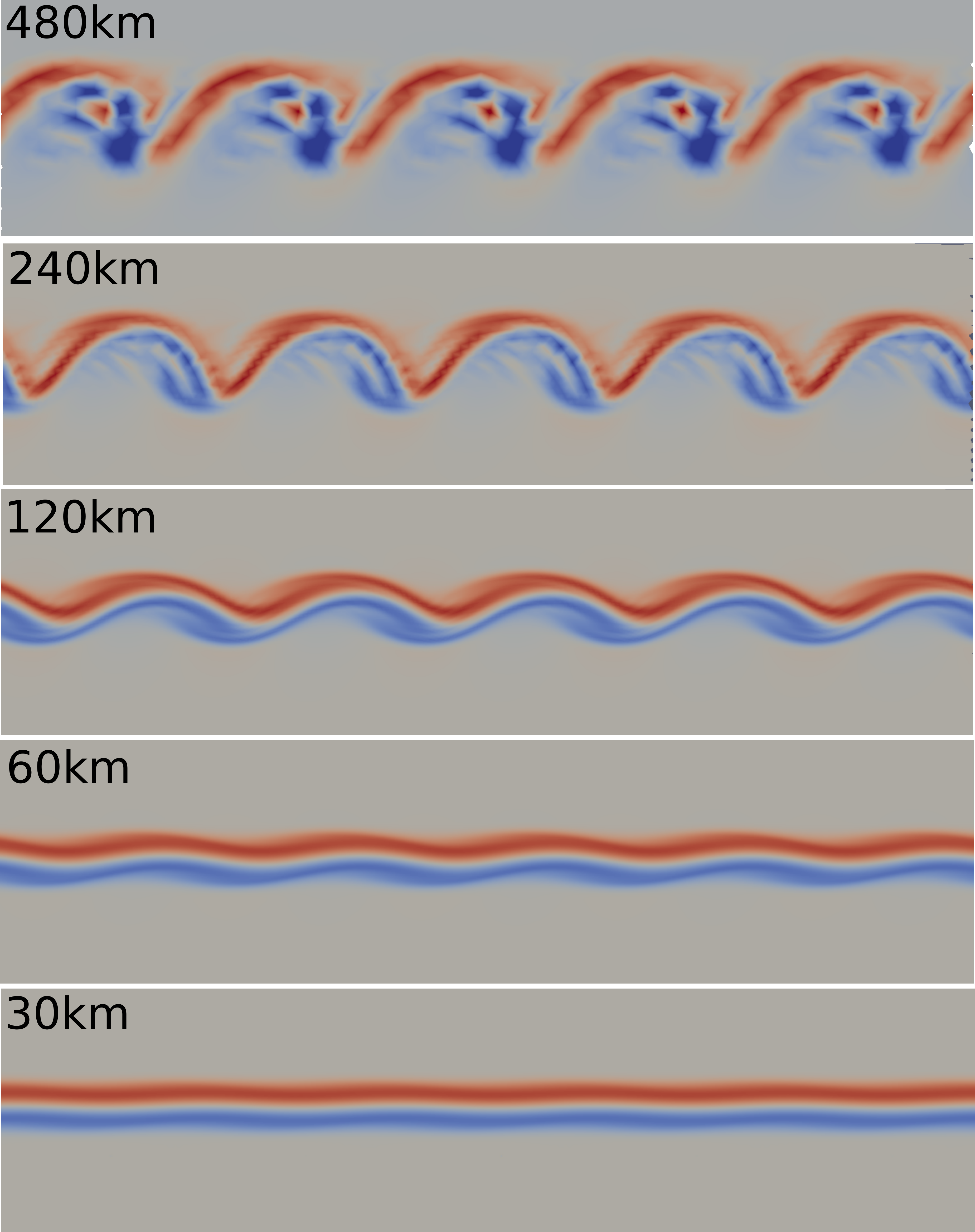}
  \caption{The vorticity field on day 5 over latitudes from $45^\circ$ to
    $90^\circ$ in the barotropic instability test case, starting from
    an unperturbed thickness field.  The 
    zonal flow structure is decently preserved by a 
    resolution of 60km and up.}
  \label{fig:tc8-unperturb}
\end{figure}

\subsection{Simulation of barotropic instabilities}
Instabilities, both barotropic and baroclinic, in the ocean and
atmosphere are responsible for the transient and diverse weather
phenomena. Therefore, in weather forecast, where accuracy is the
primary concern, it is vitally important to accurately simulate the
instabilities. But instabilities often happen erratically, if not
randomly, and accurately simulating the instabilities in the
ocean and atmosphere has been a challenge. The one-layer shallow water
model, by nature, is not capable of simulating the baroclinic
instabilities. Here, we only exam the capability of our model to
simulate the barotropic instability. We will use the familiar Galewsky
et al barotropic instability test case (\cite{Galewsky2004-hr}). This
test case starts with a perfectly balanced fast zonal jet at
$45^\circ$ north. Due to the barotropic instability, this jet will
quickly start to wobble, and eventually give rise to a stream of
eddies alongside the jet around day 5, in line with the time scale
observed in our weather system. The exact setup of this experiment can
be found in the original reference. 

\begin{figure}[h]
  \centering
  \includegraphics[width=4.6in]{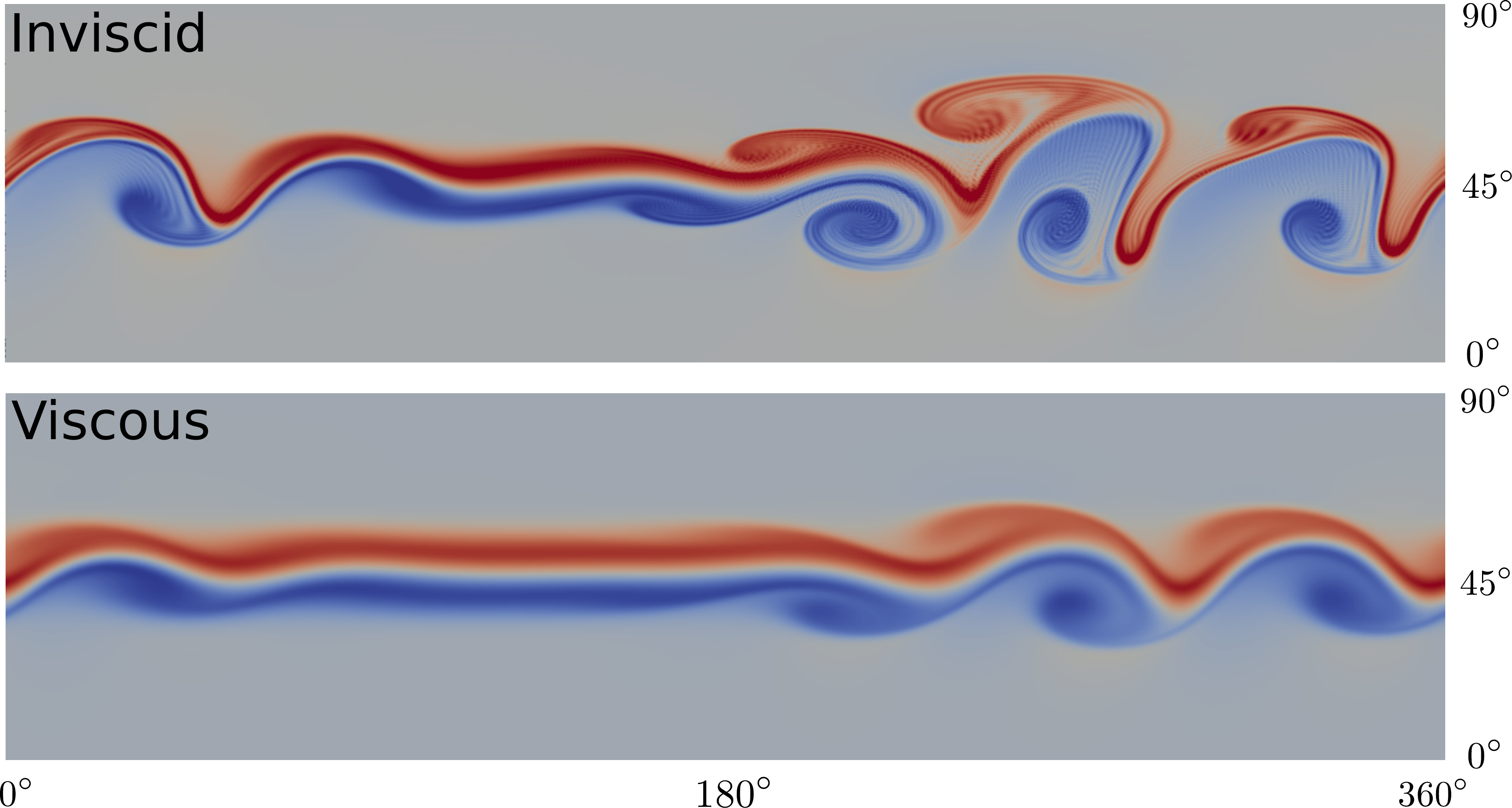}
  \caption{The vorticity fields on day 6 in the barotropic instability
    test case, with the perturbed initial thickness 
  field, on a quasi-uniform mesh with 655,362 grid cells (approx.~res.~30km). }
  \label{fig:vort-tc8}
\end{figure}

As just mentioned, 
initially, the thickness and zonal velocity field are in a perfect
geostrophic balance. A model with unlimited precision can maintain
this balance 
indefinitely. As pointed in the original article
\cite{Galewsky2004-hr}, maintaining this balance is trivial for
spectral models using spherical harmonics, because, with all the modes
initially horizontal, it will be very difficult for non-horizontal
modes to grow. This task is, on the other hand, challenging for
low-order models based on grids that are not aligned with the
latitudes, especially unstructured meshes such as ours. In these
models, the discretization errors can quickly get into the system,
and distort the geostrophic balance.
In our first
test, we run the test case using the unperturbed thickness field up to
day 5, using a sequence of quasi-uniform SCVT meshes, with mesh
resolutions ranging from 480km up to 30km. The goals are (i) to see the
effect of discretization errors at various resolutions, as measured by
the departure from the exact geostrophic balance, and (ii) to see
what resolution is required to maintain a decent geostrophic balance
up to day 5. In Figure \ref{fig:tc8-unperturb}, we plot the vorticity
field on day {5} from the simulations. In the top two panels for 480km
and 240km resolutions, the vorticity exhibit a significant depart from
the purely zonal flow, with a Rossby wave as wide as $45^\circ$ in
latitude. At 120km (third panel), which is traditionally placed in the
mesoscale range, the Rossby wave width is significantly reduced. At
60km (fourth panel), the maintenance of the geostrophic balance and
zonal flow structure is decent. 

% \begin{figure}[h]
%   \centering
%   \includegraphics[width=4.6in]{vorticity-errors-tc8.png}
%   \caption{Relative errors in vorticity (inviscid run), compared with
%     the T341 spectral solution. }
%   \label{fig:vort-err-inv}
% \end{figure}

In order to qualitatively compare the dynamics simulated by our model
with previous results, we pick the quasi-uniform $30$km mesh, as it
is shown in the above that the effect of the truncation errors on this
mesh on the overall dynamics is minimal. 
In Figure \ref{fig:vort-tc8}, we plot the vorticity field on day 6, for both
the inviscid and viscous cases, with the perturbed thickness field. In
the inviscid case, the instability is obviously more pronounced. In
both cases, the results match those of the previous studies
(\cite{Galewsky2004-hr,Chen2008-gz,Weller2012-jw}). 

Left unexplored in the work is the question of whether, in the
inviscid case, the solution actually converges to a regular solution
free of singularities (shocks). Common wisdom say that shocks are
inevitable in a compressible system, in the absence of diffusions. The
shallow water equations  is a slightly compressible system, and the
Coriolis force is known as a regularizing effect. It is not clear
whether this regularizing effect can help to suppress the emergence of
the shocks. No theoretical results exsit concerning one way or the
other. In \cite{Galewsky2004-hr}, the convergence of the solution in
the viscous case is established numerically. The issue in the inviscid
case was not explored, due to the prohibitive computational
costs. Such costs are also beyond our reach at the moment. 

\section{Discussions and conclusions}\label{sec:conclu}
%% Goal and plan for this work
This overarching goal of this  work is to evaluate the newly developed
conservative numerical scheme (EEC), so as to preliminarily assess its
potentials in real-world applications. To this end, this work focuses
on, and is organized around, a set of pre-defined properties such as
accuracy and conservation of key quantities. The evaluations make use
of a suite of old and new test cases.

%% Summary of results
% Using existing test cases for shallow water models as well as a
% specially designed test case over an arbitrary bounded domain, we find
% that the scheme is first-to-second order accurate in both the
% thickness variable and the vorticity variable, conserve both the total
% energy and potential enstrophy up to the time truncation errors. We
% also find that the scheme
% exhibits a better control of the divergence variable, sigaling that
% scheme does a better job facilitating the geostrophic adjustment
% process, and the scheme provides a more realistic representation of
% the energy and enstrophy spectra at the end of a long simulation,
% confirming that energy-and-enstrophy conserving scheme is able to
% eliminate some of the spurious energy cascade at the finest scales.

With regard to accuracy, the EEC scheme has a convergence rate between
the first and second orders, as evidenced in the SWSTC \#2 and \#5. On
bounded domains, the convergence rate is slightly lower, but remains
within that range. On the SWSTC \#2, our scheme exhibits a slightly
faster convergence rate, second {\itshape vs} first, than the extended
Z-grid scheme of Eldred (\cite{Eldred2015-yg}). This could be
attributed to the fact that we use the more accurate area-weighted
mapping throughout the scheme, while the extended Z-grid scheme uses
less accurate remapping operators for the sake of conservations. Both
our work and that of Eldred demonstrate that the accuracy of the
numerical scheme is sensitive to the mesh quality.

The EEC scheme is proven to conserve the mass and vorticity up to the
roundoff errors, and to conserve the total energy and potential
enstrophy up to the time truncation errors. These results are
confirmed by the numerical experiments conducted in this study,
including one on a bounded domain. A more intriguing question is what
advantage(s) a comprehensively conservtaive numerical scheme has in a
long-term simulations. Arakawa and Lamb (\cite{Arakawa1981-dy})
argues, via a test case with a channel flow over a ridge, argues that
a doubly conservative scheme can prevent unphysical energy cascade
towards the finest resolvable scales. Our study, via a different test
case and through a different perspective, confirms this point. We
showed that, in the SWSTC \#5, the EEC scheme, compared withan
energy-conserving only scheme (MPAS-sw), can deliver a more structured
vorticity field at the end of a long simulation (250 days), with more
physically realistic energy and enstrophy spectra.

The EEC scheme studied here has been shown to possess the optimal
dispersive wave relations among the second-order accurate numerical
schemes. It is well known that properly represented dispersive wave
relations with positive group speeds are essential to the geostrophic
adjustment process, through which large-scale geophysical flows evolve
among geostrophically balanced structures. Geostrophic balance is hard
to verify directly without the expensive reconstruction of the
velocity and pressure fields, but geostrophically balanced flows are
close to being non-divergent. Therefore, the divergence variable is a
natural choice as an indicator of the maintenance of geostrophically
balanced structures in the flows. A flow that starts out geostrophically
balanced with a small divergence variable should remain so throughout 
the simulation period. This is what is observed in the SWSTC \#5 with
the EEC scheme: the divergence variable remains small and exhibits no
discernible growth over a period of 250 days. In contrast, the
divergence variable from the MPAS-sw model, based on a C-grid scheme
with a decent representation of the dispersive wave relations, shows a
slight but unmistakenable growth over the same period. 

Effort has also been made to assess the capability of the scheme to
simulate flow dynamics close to those in the real world.
Snapshots of the vorticity field from the Galewsky et al test case
demonstrate that while the scheme has a low order of accuracy, it is
capable of simulating the barotrophic instabilities at mid- to
high-resolution, even over unstructured meshes.
%% Plan for future works
Unfortunately, the potentials of this scheme in real-world
applications cannot be 
fully assessed in its current implementation, due to the limitations
posed by the shallow water 
model, chief among them the lack of vertical variations and
convections. Hence, more work remains to be done in terms of
development and testing. On the development side, the scheme needs to
be extended and implemented for more complex models that incorporate
vertical variations and convections, such as the primitive
equations. On the testing side, real-world applications, or
well-designed test cases that are closer to the real-world
applications, need to be utitlized on the numerical scheme. These
important undertakings will be reported elsewhere.

In order to be impactful in real-world applications, the model will
have to be able to run efficiently. For this study, the model is
implemented for a single CPU, with the possible utilization of a high
performance GPU. However, we note that there is nothing in this model
that prevents or hinders a parallel implementation on heterogeneous
computing architechtures, which can be expected to multiply the
performance of the model. In the current study, the
4th order Runge-Kutta  time 
stepping scheme with a uniform time step size is used. Further
performance boost is expected from adaptive meshing and
local time stepping techniques, as demonstrated in
\cite{Hoang2019-mj}. Implementation of the scheme on modern computer
architectures, with advanced adaptive time stepping scheme, consitutes
another direction of research that will be separately pursued in
the future. 

\section*{Acknowledgements}
The first author acknowledges helpful discussions with Leo Rebholz on
the iterative scheme for the elliptic system.

% BibTeX users please use one of
%\bibliographystyle{spbasic}      % basic style, author-year citations
%\bibliographystyle{spmpsci}      % mathematics and physical sciences
\bibliographystyle{plain}      % mathematics and physical sciences
\bibliography{references}   % name your BibTeX data base

\end{document}